\documentclass[reqno]{amsart}

\usepackage{amsrefs}
\usepackage{epsf,color}
\usepackage{amsmath}
\usepackage{amscd}
\usepackage{enumerate}
\usepackage{amsfonts}
\usepackage{graphicx}
\usepackage[all]{xy}
\usepackage{mathrsfs}
\usepackage{setspace}
\DeclareFontEncoding{OT2}{}{}

\numberwithin{equation}{section}

\newcommand{\bR}{{\mathbb R}}

\newcommand{\bZ}{{\mathbb Z}}

\newcommand{\cA}{{\mathcal A}}

\newcommand{\cC}{{\mathcal C}}

\newcommand{\cI}{{\mathcal I}}

\newcommand{\cK}{{\mathcal K}}

\newcommand{\cM}{{\mathcal M}}

\newcommand{\cU}{{\mathcal U}}

\newcommand{\Fuk}{\operatorname{Fuk}}

\newcommand{\wind}{\operatorname{wd}}
\newcommand{\hol}{\operatorname{hol}}
\newcommand{\Symp}{\operatorname{Symp}}
\newcommand{\Ham}{\operatorname{Ham}}
\newcommand{\Ob}{{\mathcal Ob}}
\newcommand{\Mor}{{\mathcal Mor}}

\newtheorem{thm}{Theorem}[section]
\newtheorem{cor}[thm]{Corollary}
\newtheorem{lem}[thm]{Lemma}
\newtheorem{prop}[thm]{Proposition}
\newtheorem{defin}[thm]{Definition}

\theoremstyle{remark}
\newtheorem{rem}[thm]{Remark}

\newtheorem*{claim}{Claim}

\newcommand{\noproof}{
\begin{flushright}
\qedsymbol
\end{flushright}}

\title[On the Fukaya Categories of Higher Genus Surfaces]{On the Fukaya Categories of Higher Genus Surfaces}

\author[M.~Abouzaid]{Mohammed Abouzaid} \date{\today}

\begin{document}

\begin{abstract}
We construct the Fukaya category of a closed surface equipped with an area form using only elementary (essentially combinatorial) methods.  We also compute the Grothendieck group of its derived category.
\end{abstract}

\maketitle

\section{Introduction}
The Fukaya category  which was introduced in~\cite{fukaya1} has seen increasing interest in recent years.  Much of the interest was sparked by Kontsevich's homological mirror symmetry conjecture~\cite{kontsevich} which relates the Fukaya category of a symplectic manifold to the category of coherent sheaves of a conjectural mirror.  However, there have been recent results that underscore the relevance of this category to problems in pure symplectic topology (See ~\cites{fukaya-seidel-smith,nadler}, and ~\cite{costello} for example).

Despite this interest, there are very few compact symplectic manifolds whose Fukaya category is well understood.  Even linear invariants which can be extracted from it (of these, Hochschild (co)-homology has attracted the most interest) are only known in very few situations, and then, only as a result of a proof of homological mirror symmetry (although Picard-Lefschetz theory yields an approach which should be sufficient to understand, say, the Hochschild cohomology of the Fukaya category of Calabi-Yau complete intersections in projective space).  This means that the only compact symplectic manifolds for which computations have been performed are the elliptic curve~\cite{PZ}, ~\cite{polishchuk}, quartic surfaces~\cite{seidel-quartic}, and abelian varieties~\cite{fukaya2},~\cite{KS1}.  Also as a consequence of homological mirror symmetry, we also have a good understanding of the Fukaya-Seidel category (which can be thought of as the Fukaya category of a symplectic manifold relative a symplectic hypersurface) for some (mostly two dimensional) Lefschetz fibration~\cite{seidel-VC1},~\cite{seidel-VC2},~\cite{AKO1}, and~\cite{AKO2}.

In this paper, we consider another family of examples: Riemann surfaces.  Beyond the fact that there is a well developed theory for the study of curves on surfaces (of which we use only a small part), there are two motivations for studying these examples.  From the point of view of homological mirror symmetry, an understanding of the Fukaya category of higher genus Riemann surfaces would give at least some indication as to what form the conjecture should take in the case of complex manifolds of general type. This is a project whose main proponent for the moment is Katzarkov.  We also believe that studying this family of examples is interesting without reference to mirror symmetry, as we can hope to use them as ``building blocks'' to understand the Fukaya categories of higher dimensional symplectic manifolds.

The main result of the paper is the following
\begin{thm} \label{main}
If the Euler characteristic of $\Sigma$ is negative, then the Grothendieck group of the derived Fukaya category of $\Sigma$  is isomorphic to 
\[ H_1(S \Sigma, \bZ) \oplus \bR ,\]
where $S \Sigma$ is the unit tangent bundle.
\end{thm}
\begin{rem}
We will denote this Grothendieck group by
\[ K(\Fuk(\Sigma)) ,\]
although we should strictly speaking write $K_0(\Fuk(\Sigma))$ to emphasize the fact that we are not computing any higher $K$-groups.
\end{rem}
\begin{rem}
The map to $H_1(S \Sigma, \bZ)$ is the map which takes every smooth embedded curve $\gamma$ to the homology class of its lift to $S \Sigma$, whereas that to $\bR$ is given by computing the holonomy of a connection $1$-form, and hence is not canonical.  In reality, we will prove the theorem by choosing a (non-canonical) splitting 
\[ H_1(S \Sigma, \bZ) \cong H_1(\Sigma, \bZ) \oplus \bZ/ \chi(\Sigma)\bZ. \]
\end{rem}
\begin{rem}
The Fukaya category we shall consider will have as objects essential closed embedded curves equipped with the bounding spin structure.  We hope to return in a future paper to the subject of the Fukaya category in which $U(1)$ local systems are allowed.
\end{rem}
We end this introduction with an outline of the paper.  In Section \ref{floer-hlgy}, we define Floer homology for a certain class of immersed curves on $\Sigma$, and use this in the next section to define the Fukaya category of $\Sigma$. The next two sections prove generalizations of well known results about Hamiltonian invariance and a geometric interpretation of some extensions in term of resolving singularities.  The main technical part of the paper is Section \ref{triv-obj}, where we prove the existence of the map from $K(\Fuk(\Sigma))$ to $H_1(S \Sigma, \bZ) \oplus \bR$.  

In Section \ref{relations} and \ref{computing} we complete the proof of Theorem \ref{main} by reducing it to some well-known facts about the mapping class groups of surfaces.  The remainder of the paper consists of four appendices, in which we collect, for the convenience of the reader, some useful facts about winding numbers, the symplectomorphism groups of surfaces, the $A_{\infty}$ pre-category formalism of Kontsevich and Soibelman, and the construction of the derived category of an $A_{\infty}$ pre-category using twisted complexes.

\subsection*{Acknowledgments} 
I benefited from conversations about some aspect of this project or the other with (in alphabetical order)  Kevin Costello,  Matthew Day, Benson Farb, Ludmil Katzarkov, Andy Putman, Paul Seidel, Nathalie Wahl,  Shmuel Weinberger, and Robert Young.  I would also like to thank the Department of Mathematics of the Chinese University of Hong Kong and the Institute of Mathematical Sciences, and especially Conan Naichung Leung, for providing me with a supportive work environment during my stay in Hong Kong where part of this paper was written.  I would finally like to thank the referee for an exceptional job.

\section{Floer homology for immersed curves}\label{floer-hlgy}

Floer homology is a homology theory which assigns a group to appropriate pairs of Lagrangian submanifolds of a symplectic manifold.  The fact that these homology groups in dimension $2$ are combinatorial in nature is well known among symplectic topologists and has been used sporadically in the literature, for example in many of the aforementioned results on homological mirror symmetry for elliptic curves, or for Lefschetz fibrations in dimension two.  However, the non-expert still has to parse through the literature on $J$-holomorphic curves, with their compactness and regularity properties, in order to even get a definition of Floer homology, not to mention the proofs of basic results.

It is in order to make the paper accessible to non-experts that we therefore begin with a construction of Floer homology ab initio.  We explain enough of the theory that the relevant ideas and techniques may become clear.  We hope that those unfamiliar with Floer theory will forgive this partial state of affairs, while the expert will excuse the inclusion of detailed proofs of many results which she will think are well-known and established.  Maybe she will accept the excuse that Floer theory usually restricts itself to embedded Lagrangians, whereas we will need to use some immersed curves in our arguments.  It is possible to replace the use of immersed curves by more algebraic arguments, but that would make the details more cumbersome.

Let $\Sigma$ be a surface, possibly with boundary, equipped with a symplectic form $\omega$.

\begin{defin}
An {\bf unobstructed} curve $\gamma \subset \Sigma$ is a properly immersed oriented smooth curve which is the image of a properly embedded line in the universal cover of $\Sigma$. 
\end{defin}

As we will see later, our curves will always be equipped with one marked point.  Strictly speaking, we only care about closed curves.  However, as we will be passing to infinite covers in some of our arguments, we will also need to consider the non-closed case.  For convenience, we will adopt notation from homological algebra, and write $\gamma[1]$ for the curve obtained by changing the orientation of $\gamma$.
\begin{figure}
\begin{center}
$\begin{array}{c} \epsfxsize=1in
\epsffile{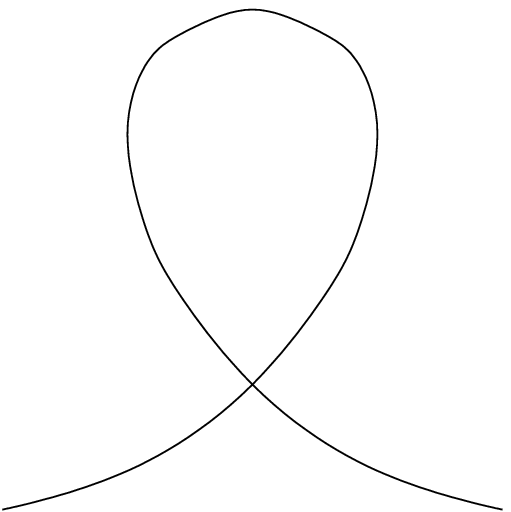}
\end{array}$
\end{center}
\caption{A ``fish-tail''}
\label{fish}
\end{figure}
\begin{lem}
A properly immersed smooth curve $\gamma$ is unobstructed if and only if the following conditions hold:

\begin{enumerate}
\item $\gamma$ is not null-homotopic, and
\item $\gamma$ does not bound an immersed ``fish-tail'' (See Figure \ref{fish}).
\end{enumerate}
\end{lem}
\begin{rem}
Formally, a ``fish-tail'' is a map from the disc $D$ to the surface $\Sigma$ which takes the boundary of $D$ to $\gamma$ and which is non-singular at every point except one point on the boundary. \end{rem}
\begin{proof}
Note that $\gamma$ is null-homotopic if and only if it is the image of an immersed circle in the universal cover, and hence cannot be unobstructed.  If $\gamma$ is not null homotopic, then it lifts to a properly immersed curve $\tilde{\gamma}$ in the universal cover.  The curve $\tilde{\gamma}$ is embedded if and only if it does not bound a fishtail, which would project to give an immersed fishtail in $\Sigma$.  For the converse, choosing the basepoint to be the singular point of the fishtail, it is easy to see that an immersed fishtail in $\Sigma$ with boundary on $\gamma$ lifts to such a fishtail in the universal cover with boundary on $\tilde{\gamma}$.
\end{proof}

Let $\gamma_1$ and $\gamma_2$ be unobstructed curves which intersect transversally.  Let $c$ be such a transverse intersection point.  We will think of $c$ as a ``morphism from $\gamma_1$ to $\gamma_2$.''  We always have an orientation preserving map from a neighbourhood of $p$ to a neighbourhood of the origin in $\bR^{2}$ taking $\gamma_1$ to the $x$-axis and $\gamma_2$ to the $y$-axis.  We call such a map a {\bf local model} for the intersection $p$.

\begin{defin}
The {\bf degree} of $c$ is $1$ if there is a local model for which the natural orientations of the $x$ and $y$ axes and the orientations of $\gamma_1$ and $\gamma_2$ agree.  The degree of $c$ is $0$ otherwise (See Figure \ref{intersect-degree}).
\end{defin}

\begin{rem}
Note that the degree of $c$ as a  ``morphism from $\gamma_1$ to $\gamma_2$'' is always different from its degree as a  ``morphism from $\gamma_2$ to $\gamma_1$.''  We are therefore abusing notation by assigning a degree to $c$.  In practice, however, this should cause no confusion as the intended meaning will be clear from the context.
\end{rem}
\begin{figure}
\begin{center}
$\begin{array}{c@{\hspace{.5in}}c} \input{degree_one.pstex_t} &
\input{degree_zero.pstex_t} \\
\mbox{Degree  1 } & \mbox{Degree 0}
\end{array}$
\end{center}
\caption{}
\label{intersect-degree}
\end{figure}

Ideally, we would define the Floer complex to be a complex over $\bR$ generated by the intersection points. We would be able to define the differential since it involves only a finite sum of contributions given by certain bigons.  However, the higher products which we will discuss in Section \ref{infty_structure} may involve contributions of countably many polygons, and it is not clear that the corresponding series converge.  It is in order to address these convergence problems that the use of Novikov rings was introduced in symplectic topology.
\begin{defin}
The {\bf universal Novikov ring} $\Lambda$ is the ring of formal power series of the form
\[ \sum_{i \in \bZ} a_{i} t^{\lambda_i} \]
where the coefficients $a_{i} \in \bZ$ vanish for all sufficiently negative $i$, and $\lambda_i$ are a sequence of real numbers satisfying \[\lim_{i \to \infty} \lambda_i = \infty . \]
\end{defin}

The idea is that the coefficients $\lambda_i$ keep track of areas of polygons, and that if there are infinitely many polygons that appear in some computation, there are still only finitely many whose area is below a certain bound.  In particular,  the coefficients $a_i$ correspond to the finite number of polygons of area $\lambda_i$.  We now begin implementing these ideas:

\begin{defin}
The {\bf Floer complex} of $(\gamma_1, \gamma_2)$ is the free $\bZ / 2 \bZ$-graded module over the universal Novikov ring $\Lambda$ generated by intersection points of $\gamma_1$ and $\gamma_2$. It is denoted
\[ CF^{i}(\gamma_1, \gamma_2) = \bigoplus_{\deg(c) = i} [c] \cdot \Lambda .\]
\end{defin}

 We define a differential 
\[ d: CF^{i}(\gamma_1, \gamma_2)  \to CF^{i+1}(\gamma_1, \gamma_2) \]
by counting immersed bi-gons.  Formally, we have
\begin{defin}
Let $p$ and $q$ be intersection points between $\gamma_1$ and $\gamma_2$.  We define
\[ \cM(q,p) \]
to be the set of immersed discs $(u, \partial u)$ which satisfy the following conditions (See Figure \ref{differential}):
\begin{itemize}
\item The points $p$ and $q$ lie on $\partial u$, and
\item If we give $(u, \partial u)$ the induced orientation from the inclusion into $\Sigma$, the clockwise path from $p$ to $q$ along $\partial u$ is an immersed path lying on $\gamma_1$, while the counterclockwise path lies on $\gamma_2$, and
\item Both $p$ and $q$ are convex corners of the disc $u$.
\end{itemize}
\end{defin}

Strictly speaking, we should think of $u$ as an orientation preserving map
\[ u: (D, \partial D) \to (\Sigma, \gamma_1 \cup \gamma_2) \]
which maps $1$ to $p$, $-1$ to $q$ and satisfies the above conditions. In particular, we can assume that $u$ is an immersion away from $\pm 1$.  Note that we would then have to consider such maps only up to smooth reparametrization of the disc.  We will sometimes use this point of view for additional clarity.

\begin{rem} \label{quadrants}
If we choose a local model near $q$ taking $\gamma_1$ to the $x$-axis, and $\gamma_2$ to the $y$-axis, then for every $v \in \cM(q,p)$  the image of a neighbourhood of $-1$ is contained in either in the first quadrant ($x<0$ and $y>0$) or the third quadrant ($x>0$ and $y<0$).  Similarly, every $ u \in \cM(r,q)$ must locally lie entirely in one of the remaining quadrants in a local model at $q$.
\end{rem}
%The following lemma is intuitively obvious
%\begin{lem}\label{just_two}
%An immersed disc is determined, if it exists, by which quadrant contains it in a local model at both of its corners. In particular, $\cM(q,p)$ consists of at most two elements.
%\end{lem}
%\begin{proof}
%Choosing $p$ as a basepoint of $\Sigma$, let us pass to the universal cover where the lifts $\tilde{\gamma}_1$ and $\tilde{\gamma}_2$ are embedded.  These lifts intersect at $\tilde{p}$.  Let us fix quadrants at $\tilde{p}$ and at $\tilde{q}$, and consider the intersection of the half rays originating at $\tilde{p}$ and $\tilde{q}$ along  $\tilde{\gamma}_1$ and $\tilde{\gamma}_2$.  If the intersections are both non-empty, they form a circle in $\bR^2$.  The immersion condition requires that every element of $\cM(p,q)$ has this circle as its boundary, but it is easy to see that an immersed circle in $\bR^2$ bounds at most one immersed disc (up to reparametrization).
%\end{proof}

\begin{figure}
\begin{center}
$\begin{array}{c} \input{differential2.pstex_t}
\end{array}$
\end{center}
\caption{An element of $\cM(q,p)$.}
\label{differential}
\end{figure}

We can now define the differential according to the formula
\[ d [p] = \sum_{q, u \in \cM(q,p)} (-1)^{s(u)} t^{\omega(u)} [q] .\]

Where we write $\omega(u)$ for the area of $u$ given the induced orientation as a subset of $\Sigma$.  In particular, $\omega(u)$ is always positive.  The sign $(-1)^{s(u)}$ is determined as follows:  The sign is the product of a $-1$ contribution coming from every marked point which appears on the boundary of $u$, and $-1$ contribution whenever the given orientation of $\gamma_2$ disagrees with its orientation as the boundary of $u$.

%By Lemma \ref{just_two} the coefficients in the above sum never have absolute value greater than $1$.

\begin{rem}
It is fairly easy to see that if $\cM(q,p)$ is non-empty, then the degree of $p$ and $q$ differ by $1$.  One can also see that the differential may act non-trivially on both $ CF^{0}(\gamma_1, \gamma_2)$ and  $ CF^{1}(\gamma_1, \gamma_2)$.  
\end{rem}

The following lemma justifies our non-standard definition of unobstructed curves:
\begin{lem} \label{square_zero_proof}
The Floer differential satisfies $d^{2} = 0$.
\end{lem}
\begin{proof}
In Floer theory, this is proved by considering discs which occur in $1$-dimensional moduli.  In our case, this means considering immersed discs with one non-convex corner.

Let us choose a generator $[p]$.  We have:
\[ d^{2}[p] = \sum_{q, u \in \cM(q,p)} \sum_{r, v \in \cM(r,q)} (-1)^{s(u)+s(v)} t^{\omega(u) + \omega(v)} [r] .\]

We must show the vanishing of the coefficient of every generator $[r]$ in the right hand side of the above equation.  In particular, we will show that the above coefficients cancel pairwise.

Let us choose an intersection point $q$ such that there exists $u \in \cM(q,p) $ and $v \in  \cM(r,q) $.  By Remark \ref{quadrants}, there are four possibilities for which quadrants $u$ and $v$ occupy near $q$.  The argument is the same in all cases, so let us assume that $u$ lies in the first quadrant, and $v$ in the fourth.

Since $\partial u$ and $\partial v$ both contain segments from $\gamma_1$, and these segments overlap near $q$, we know that one of these segments must be contained in the other.  Again, there are two cases here, and we will simply address the case where $r$ lies on the boundary of $u$.  In this case, the tangent vector of $\gamma_2$ at $r$ must point towards the interior of $u$.  This forces $\gamma_2$ to intersect the boundary of $u$ transversally at yet another point (See Figure \ref{square_zero}).  However, this point cannot lie on the part of $\partial u$ which lies on $\gamma_2$ since this would create an immersed fish-tail.  We must therefore have yet another intersection of $\gamma_2$ with $\gamma_1$ along the boundary of $u$.  Let us denote this additional intersection by $q'$.

We now have two elements $u' \in \cM(q',p)$ and $v' \in \cM(r,q')$.  It is clear that $\omega(u) + \omega(v) = \omega(u') + \omega(v')$.  Further, if there are no marked points on the boundary of these four curves, $s(u)$ and $s(u')$ have the same parity, while the parities of $s(v)$ and $s(v')$ must differ, so that the total contribution of these four discs to the coefficient of $r$ in $d^2[p]$ vanishes.  Note that any marked point contributes to two discs, so the result holds even in their presence.
\end{proof}

\begin{figure}
\begin{center}
$\begin{array}{c@{\hspace{.5in}}c} \input{square_zero_1.pstex_t} &
\input{square_zero_2.pstex_t}
\end{array}$
\end{center}
\caption{}
\label{square_zero}
\end{figure}

\begin{rem}
In order to avoid the case by case argument that we used, we could observe that when we let $w$ be the union of $u$ and $v$, we obtain a new immersed polygon with boundary on $\gamma_1 \cup \gamma_2$, with corners at $p$ and $r$.  However, one of these corners is not be convex.  Say the non-convex corner is $p$.  Then there must be two segment starting at $p$ which intersect $\partial w$ transversally, one along $\gamma_1$ and the other along $\gamma_2$.  One of these segments recovers the decomposition of $w$ as $u \cup v$, whereas the other yields a new decomposition of $w$ as two polygons $u'$ and $v'$.  These correspond to two terms whose contributions to $d^2[p]$ cancel each other. 
\end{rem}
\begin{defin}
The {\bf Floer cohomology} of $\gamma_1$ and $\gamma_2$ is the cohomology of the complex $CF^{*}(\gamma_1, \gamma_2)$ with respect to the differential $d$.  It is denoted
\[ HF^{*}(\gamma_1, \gamma_2) .\]
\end{defin}

\section{$A_{\infty}$ structure} \label{infty_structure}
The next step in constructing the Fukaya category is to describe how to compose morphisms.  However, such a composition will only be associative up to homotopy.  There is an algebraic theory of such categories which are only associative up to homotopy, but the theory works best whenever there is in fact a whole sequence of higher homotopies which correct the failure of associativity of the product.
Let us therefore consider $\{ \gamma_{i} \}_{i=0}^{k}$,  unobstructed curves in general position (It is enough to assume that there are no triple intersections).  

\begin{figure}
\begin{center}
 \input{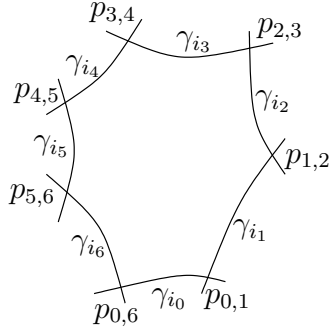} 
\end{center}
\caption{A polygon in $\cM( p_{0,6}, p_{5,6}, \cdots,p_{0,1})$}
\label{definition_polygon}
\end{figure}

We define a product 
\[ m_{k}:  CF^{*}(\gamma_{k-1}, \gamma_{k}) \otimes \cdots \otimes  CF^{*}(\gamma_{1}, \gamma_{2}) \otimes CF^{*}(\gamma_{0}, \gamma_{1})  \to CF^{*}(\gamma_{{0}}, \gamma_{k}) \]
by counting immersed $k+1$-gons all whose corners are convex.  More precisely, let $p_{i,i+1}$ be an intersection point between $\gamma_{i}$ and $\gamma_{i+1}$.
\begin{defin}  Denote by
\[ \cM( p_{0,k}, p_{k-1,k}, \cdots,p_{0,1}) \] the set of immersed polygons $u$ such that the boundary of $u$ is a union of segments 
\[ [p_{i-1,i}, p_{i, i+1}] \subset \gamma_i \]
with convex corners at every $p_{i,i+1}$ and such that the natural orientation of the boundary of $u$ corresponds to increasing $i$ (See Figure \ref{definition_polygon}).
\end{defin}
\begin{rem}
The orientation condition dictates that in a local model near each $p_{i,i+1}$ taking $\gamma_i$ to the $x$-axis, and $\gamma_{i+1}$ to the $y$-axis, such a polygon lies in the second or fourth quadrant.  At $p_{0,k}$, we must consider the opposite local model in order for $u$ to be contained in the even quadrants.  %However, the analogue of Lemma \ref{just_two}, is false in this situation.
\end{rem}

 The contribution of such polygons is given by
\begin{equation}\label{polygons_infty} m_{k} ( [p_{k-1,k}], \cdots,  [p_{0,1}]) = \sum_{p_{0,k}} \sum_{ u \in  \cM(p_{0,k},p_{k-1,k}, \cdots , p_{0,1})} (-1)^{s(u)} t^{\omega(u)} [p_{0,k}]  . \end{equation}

Where $\omega(u)$ is the area of $u$ as before, and the sign $s(u)$ is determined as follows:  Assign to every odd degree intersection $p_{i,j}$ where $i<j$ a sign $(-1)^{s(p_{i,j})}$.  This sign will be positive if and only if the orientation of $\gamma_j$ as the boundary of $u$ agrees with the given orientation of $\gamma_j$.  The total sign $(-1)^{s(u)}$ is the product of the signs assigned to all odd degree intersections and to all marked points which appear on the boundary of $u$.

\begin{lem} \label{degree_product}
The degree of $m_k$ is equal to $k\mod(2)$.
\end{lem}
\begin{proof}
Since $u$ lies in one of the even quadrants, agreement of the orientation of $\partial u$ with $\gamma_i[a]$ implies agreement of the orientation of $\partial u$ with $\gamma_{i+1}[a + 1 + \deg(p_{i,i+1})]$. Keeping in mind the different convention for the last point, going all the way around the boundary of the polygon implies
\[ a_0 = a_0 + \deg(p_{0,k}) + \sum_{0 \leq i < k} 1 + \deg(p_{i,i+1}) \mod(2) ,\]
which proves
\[ \deg(p_{0,k}) =  k + \sum_{0 \leq i < k} \deg(p_{i,i+1}) \mod(2) .\]
\end{proof}

When $k=1$, we simply recover the differential which we constructed in the previous section.  The only non-trivial part is that the signs agree, but if $\cM(q,p)$ is non-empty, then either $p$ or $q$ (but not both) has even degree, hence contributes a sign to $d[p]$ which is exactly the sign which we used in the definition of the differential.  When $k=2$, we obtain the product which we claimed would only be associative up to homotopy.  First, we explain why the product respects the differential.  In other words, we prove:

\begin{lem}
If $\gamma_0$, $\gamma_1$, and $\gamma_2$ are transverse unobstructed curves, then
\[ m_2:  CF^{*} (\gamma_1, \gamma_2)  \otimes CF^{*} (\gamma_0, \gamma_1) \to CF^{*}(\gamma_0, \gamma_2) \]
satisfies the formula:
\[ m_1(m_2([q],[p])) = (-1)^{\deg(p)} m_2(m_1([q]), [p]) - m_2([q], m_1([p])) .\] 
\end{lem}
\begin{rem}
The reader might have expected the usual conventions for differential graded algebras for which the sign in the second term of the right hand-side is $(-1)^{\deg(q)}$.  One can recover the usual sign conventions by defining a new differential $\partial (q) = (-1)^{\deg(q)} m_1(q)$ and a new product $q \cdot p = (-1)^{\deg(p)} m_2(q,p)$.
\end{rem}
\begin{proof}
We rewrite the above formula as
\[  m_1(m_2([q],[p])) - (-1)^{\deg(p)} m_2(m_1([q]), [p]) +  m_2([q], m_1([p])) =0 ,\]
and expand the left hand-side into the sum of three series
\begin{align*}
\sum_{r,u \in \cM(r,q,p)} \sum_{s,v \in \cM(s,r)} (-1)^{s(u) + s(v)} t^{\omega(u) + \omega(v)} [s] \\ + (-1)^{\deg(p)} \sum_{r,u \in \cM(r,q)} \sum_{s,v \in \cM(s,r,p)} (-1)^{s(u) + s(v)} t^{\omega(u) + \omega(v)} [s] \\
  + \sum_{r,u \in \cM(r,p)} \sum_{s,v \in \cM(s,q,r)} (-1)^{s(u) + s(v)} t^{\omega(u) + \omega(v)} [s] \end{align*}
We must now show that the total coefficient of $s$ in the above expression vanishes.  As in the proof of Lemma \ref{square_zero_proof} one in fact shows that the terms cancel pairwise. As in the remark following that lemma, we should therefore consider polygons with one non-convex corner. The only non-trivial part is to check that the signs of the two terms that such polygons contribute are always opposite.  We will refer to the proof of Lemma \ref{infty_relation} for the relevant argument.
\end{proof}

In order to prove that the product is associative up to homotopy, we must show that there exists a map 
\[ h: CF^{*}(\gamma_2, \gamma_3) \otimes  CF^{*}(\gamma_1, \gamma_2) \otimes  CF^{*}(\gamma_0, \gamma_1) \to   CF^{*}(\gamma_0, \gamma_3) \]
such that
\begin{align*} m_2([r],m_2([q],[p])) - (-1)^{\deg(p)} m_2(m_2([r],[q]),[p]) = 
\\ - m_1(h([r],[q],[p])) - (-1)^{\deg(q)+\deg(p)} h(m_1([r]), [q],[p]) 
\\+ (-1)^{\deg(p)} h([r], m_1([q]),[p]) -   h([r], [q], m_1([p])) .\end{align*}
Note that the left hand-side would vanish if the product was strictly (graded) associative, and that the right hand-side should be interpreted as applying the total differential to $h$.  In fact, we have already constructed $h$; it is equal to $m_3$.  A generalization of this result is called the  $A_{\infty}$ relation.
\begin{lem} \label{infty_relation}  The higher products $m_k$ satisfy the equation
\[ \sum_{l , d, j} (-1)^{\maltese_{d}} m_{l} ( [p_{k-1, k}]  , \cdots, m_{j} ([p_{d+j, d+j+1}], \cdots , [p_{d, d+1}] ), \cdots ,[p_{0,1}]) = 0 ,\]
where
\[ \maltese_{d} =  d+\sum_{i=1}^{d} \deg(p_{i-1,i}). \]
\end{lem}
\begin{proof}
The proof is essentially the same as the proof that $m_1^2 = 0$, or that $m_2$ respects the differential.  It will suffice to explain why, given a polygon with one non-convex corner, the two contributions of the decomposition into a pair of convex polygons cancel each other.  As usual, the only issue is that of signs.

Let $w$ be such a polygon with decompositions into $u \cup v$ and $u' \cup v'$ as in Figure \ref{infty_relation_picture} with additional corners $p_{d,d'+1}$ and $p_{d',d''}$.  We must show that
\[ d+\sum_{i=1}^{d} \deg(p_{i-1,i}) + s(u) + s(v) + d' +\sum_{i=1}^{d'} \deg(p_{i-1,i}) + s(u') + s(v') = 1 \mod(2).\]

\begin{figure}
\begin{center}
 \input{infty_relation2.pstex_t}
\end{center}
\caption{}
\label{infty_relation_picture}
\end{figure}

Recall that $s(u)$ was defined as the sum of the contributions of the corners of $u$ which have degree equal to $1$ multiplied by the contribution of all the marked points.  Since all the marked points appear an even number of times, they can be safely ignored, so we assume that none exist.  In particular, every corner other than $p_{d,d'+1}$, $p_{d',d''}$, and $p_{d',d'+1}$ have the same contribution to $s(u) + s(v)$ and $s(u') + s(v')$.  From Figure \ref{induced_orientations}, we see that since the outgoing edge at $p_{d,d'+1}$ appears with opposite orientations in the discs appearing in the two formulae, this intersection points contributes its degree to $s(u) + s(v)$.  Similarly, $p_{d',d'+1}$ contributes its degree to $s(v) + s(v')$.  On the other hand, the outgoing edge at $p_{d',d''}$ appears twice with the same orientation, so its total contribution to $s(u') + s(v')$ vanishes.  Therefore, it suffices to show that

\[ d+\sum_{1 = i}^{d} \deg(p_{i-1,i}) + \deg(p_{d,d'+1}) + \deg(p_{d',d'+1}) + d' +\sum_{i=1}^{d'} \deg(p_{i-1,i})  = 1 \mod(2).\]

But this follows easily from Lemma \ref{degree_product} which implies that
\[ \deg(p_{d,d'+1}) = (d'+1-d) + \sum_{ d < i \leq d'+1}  \deg(p_{i-1,i}) .\]
\begin{figure}
\begin{center}
$\begin{array}{c@{\hspace{.25in}}c} 

\input{induced_orientation_v.pstex_t} & \input{induced_orientation_v_prime.pstex_t}\\
\input{induced_orientation_u.pstex_t} & \input{induced_orientation_u_prime.pstex_t}
\end{array}$
\end{center}
\caption{}
\label{induced_orientations}
\end{figure}

\end{proof}

It would seem that that we have constructed an $A_{\infty}$ category (See Appendix \ref{precategories} for a definition). However, the standing assumption in this section has been that the curves $\gamma_i$ are transverse.  While this is a generic assumption, it certainly does not hold for all curves, so we have not even defined the space of endomorphisms of an object leaving us in a situation where it does not make sense to speak of units.

There are several solutions to this problem.  The most elegant and geometric solution is to introduce a ``Morse-Bott'' version of the theory, which is done for Lagrangian Floer theory in \cite{FOOO}.  This theory is in fact of manageable complexity in dimension $2$, but developing it with full rigour is still a non-trivial task.  Instead, we will opt for hiding the problem in the algebraic structure by working with the $A_{\infty}$-pre-categories of Kontsevich and Soibelman \cite{KS1}.  We include a brief review in Appendix \ref{precategories}.  The reader should be warned that the simplicity of working with $A_{\infty}$ pre-categories is deceiving, and that we will pay the price in Section \ref{triv-obj} for having chosen not to develop a ``Morse-Bott'' theory of combinatorial Floer homology.  In particular, the proof of Lemma \ref{bounding_chain} is essentially a definition in the latter theory.

\begin{defin}
Let $\Sigma$ be a surface equipped with a symplectic form.  The {\bf Fukaya category } of $\Sigma$ is the $A_{\infty}$-pre-category $\Fuk(\Sigma)$ given by the following data:
\begin{itemize} \label{define-fuk}
\item The objects are all embedded unobstructed curves in $\Sigma$.
\item The transverse subsets $\Ob(\Fuk(\Sigma))^{tr}_{n}$ are given by $n$-tuples of curves such that no three curves have a common point.
\item Among the transversal objects, the higher compositions are given as in Equation \eqref{polygons_infty}.
\end{itemize}
\end{defin}

We remark that unlike the case of a general $A_{\infty}$ pre-category, a reordering of transverse sequence is still transverse, so we may simply speak of transverse sets.  The rest of the paper is devoted to developing geometric techniques to compute the $K$-group (See Appendix \ref{twisted-and-K}) of the Fukaya category of a surface.

\section{Quasi-isomorphic curves}
Let us choose a connection $1$-form $A$ as in Appendix \ref{wind_holo}.  In this section, we prove the following proposition, which is the appropriate generalization, in the presence of immersed curves, of the Hamiltonian invariance of Floer homology: 
\begin{prop}\label{invariance}
Assume that $\gamma_1$ and $\gamma_2$ are transverse unobstructed curves which are isotopic and satisfy
\[\hol_A(\gamma_1) = \hol_A(\gamma_2) .\]
Then $\gamma_1$ and $\gamma_2$ are quasi-isomorphic objects of $\Fuk(\Sigma)$.
\end{prop}

We first prove the result assuming that $\Sigma$ is the cylinder, that $\gamma_1$ and $\gamma_2$ are circles representing the generator of the fundamental group.  Note that in this case, $\gamma_1$ and $\gamma_2$ are embedded.
\begin{figure}
\begin{center}
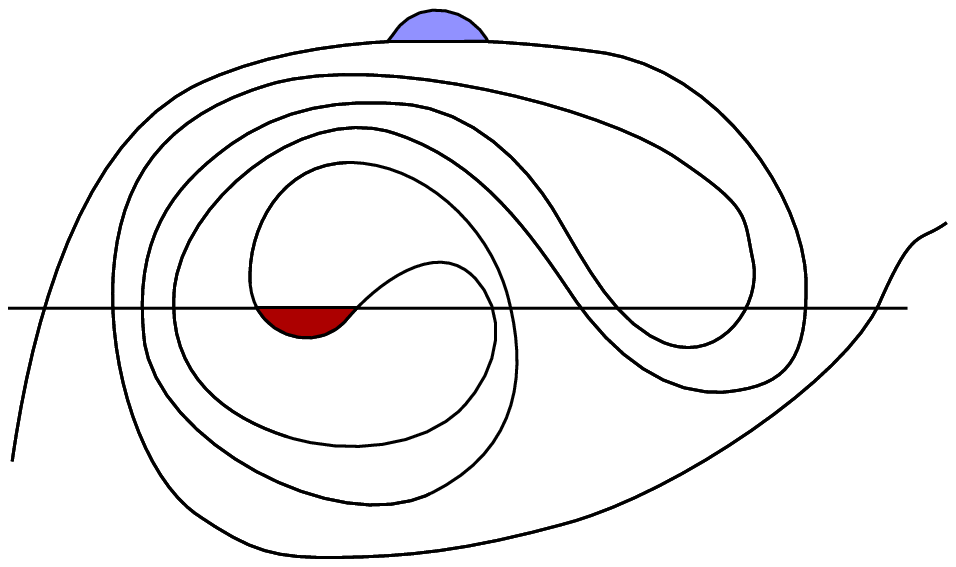 
\end{center}
\caption{}

\label{intersection-reduction}
\end{figure}

\begin{lem}
If $|\gamma_1 \cap \gamma_2|>2 $ then there exists a curve ${\gamma^{1}_2}$ which intersects $\gamma_2$ in two points, intersects $\gamma_1$ in $|\gamma_1 \cap \gamma_2|-2$ points and has the same holonomy as $\gamma_2$.
\end{lem}
\begin{proof}
Since $\gamma_1$ and $\gamma_2$ represent the same homology class, their algebraic intersection number vanishes. Consider a pair $(p,q)$ of intersection points which are the corners of an embedded disc $u$ of minimal area $a$. One element of the pair will be a positive intersection point; the other a negative intersection.  We will construct the curve $\gamma^{1}_2$ as a small perturbation of a piecewise-smooth curve $\gamma$ which agrees with $\gamma_2$ almost everywhere, except near $u$ where it will agree with $\gamma_1$ and near a yet unspecified point where it will bound a disc $u'$ of area $a$ along with $\gamma_2$ in such a way that the holonomy along $\gamma$ agrees with the holonomy along $\gamma_2$; see Figure \ref{intersection-reduction}.

We must choose $u'$ in such a way as to ensure that $\gamma'_2$ does not pick up any new intersection with $\gamma_1$ and $\gamma_2$.  Consider a trivialization of the cylinder as $S^1 \times \bR$ such that $\gamma_1$ is the fiber over $0$ of the projection map to $\bR$, and let $(m,M)$ be respectively points on $\gamma_2$ which have minimal and maximal $\bR$ values.  Note that these are points where $\gamma_2$  is tangent to a translate of $\gamma_1$.
\begin{claim}
The orientation of $\gamma_2$ at $m$ (and $M$) agrees with the translate of $\gamma_1$.
\end{claim}
\begin{proof}
By lifting to $\bR^2$ with basepoint $\tilde{m}$ over $m$, we obtain a line $\tilde{\gamma_2}$  which is tangent to some $y=c$ line at $\tilde{m}$.  If the orientation of $\gamma_2$ and the translate of $\gamma_1$ do not agree at $m$, the curve $\tilde{\gamma_2}$ ``comes from'' $x=-\infty$, has tangent vector $ - \frac{\partial}{\partial x}$ at $\tilde{m}$, then heads towards $x= + \infty$ without crossing itself (no self intersections) or the $y=c$ lines (since $c$ was the minimal $\bR$ value).  This is impossible.
\end{proof}

We can now modify $\gamma_2$ near $m$ (or $M$ depending on the relative orientations of $\gamma_2$ and the boundary of $u$) to produce a curve $\gamma$ which has the same holonomy as $\gamma_2$ (see Figure \ref{intersection-reduction}).  The desired conditions about intersection points follow from the fact that $u$ was the disc of minimal area.
\end{proof}
\begin{cor} \label{inductive_reduction}
There exists a sequence $\{\gamma_2= \gamma^{0}_2, \ldots,  \gamma^{n}_2 = \gamma_1 \} $ of curves having the same holonomy, such that $|  \gamma^{i}_2 \cap  \gamma^{i+1}_2 | =2$.
\end{cor}

In particular, since quasi-isomorphisms are transitive it will suffice to prove that successive curves in the above sequence are quasi-isomorphic.  Since it will be needed later, we will be constructing the explicit quasi-isomorphism in a slightly more general context:

Assume all the immersed discs bounded by $\gamma_1$ and $\gamma_2$ are in fact embedded with disjoint interiors;  this finite set of bigons $\{ u_k \}_{k=0}^{2n-1}$ satisfies
\[ \sum_{k} (-1)^{k} \omega(u_k) = 0 .\]
The bigons $u_k$ are numbered such that $k$ increases as one traverses the boundary of $\gamma_1$ according to its orientation. Let $\{ p_i \}_{i=0}^{n-1}$ be the corners of $u_k$ which have degree $0$ in $CF^{*}(\gamma_1, \gamma_2)$ and let $q_i$ be those corners of degree $1$.  Since the position of the two marked points does not affect the quasi-isomorphism type we may assume that they both lie on the boundary of the same disc $u_k$ so that their signed contributions cancel.  Working $\mod(2n)$, we have:
\[ m_1(p_i) = - t^{\omega(u_{2i-1})} q_{i-1} + t^{\omega(u_{2i})} q_{i} .\]
In particular, we compute:
\begin{align*}
m_1( \sum_{i=0}^{n-1} t^{ \sum_{k=0}^{2i-1} (-1)^{k} \omega(u_k) } p_i ) & = \sum_{i} - t^{ \sum_{k=0}^{2i-2} (-1)^{k} \omega(u_k)} q_{i-1} +  t^{ \sum_{k=0}^{2i} (-1)^{k}  \omega(u_k)} q_{i} \\
& = - t^{\omega(u_{2n-1})} q_{n-1} + t^{ \sum_{k=0}^{2n-2} (-1)^{k} \omega(u_k)} q_{n-1} \\
& = 0.
\end{align*}
This proves that
\[ e_{1,2} =  \sum_{i=0}^{n-1} t^{ \sum_{k=0}^{2i-1} (-1)^{k} \omega(u_k) } p_i \] is a closed morphism.  In the rest of this section, we will prove that it is a quasi-isomorphism.  First we will choose a curve $\gamma'_1$ which is $C^{1}$-close to $\gamma_1$ and which satisfies $\hol_{A}(\gamma'_1) = \hol_{A}(\gamma_1)$.  We assume that all bigons bounded by $\gamma_1$ and $\gamma'_1$ have the same area, and construct morphisms $e_{1,1'}$ and $e_{1',1}$ which we will prove are quasi-isomorphisms.  We begin with

\begin{lem} \label{strict_units}
If $\{ \alpha_n \} \cup \gamma_1$ is a finite collection of mutually transverse curves and $\gamma'_1$ satisfies the above conditions then
\[ m_k([x_k], \cdots, [x_{d+2}], e_{1,1'}, [x_{d}], \cdots [x_1]) = 0 \]
whenever $k >2$ and $x_i \in \alpha_{n_i} \cap \alpha_{n_{i+1}}$.
\end{lem}
 \begin{rem}
Of course, the $C^{1}$-closeness of $\gamma_1$ and $\gamma'_1$ is with respect to the curves $\alpha_n$.
\end{rem}
\begin{proof}
We will in fact prove that the sets
\[ \cM(x',x_k, \ldots, x_{d+2}, p, x_d, \ldots, x_1) \]
are empty whenever $k > 2$ and $p$ is an intersection point between $\gamma_1$ and $\gamma'_1$ which yields a degree $0$ morphism from $\gamma_1$ to $\gamma'_1$.  In order to do this, we begin by observing that an immersed polygon $u$ in one of the above sets agrees, near $p$ with one of the bigons $w$ which $\gamma_1$ and $\gamma'_1$ bound.  Assume that $q$ is the other corner of $w$.  If the segment along $\gamma_1$ from $p$ to $x_d$ and the segment along $\gamma'_1$ both pass through $q$, then we see that $u$ must in fact agree with the bigon from $p$ to $q$, since otherwise the inverse image of $q$ would contain a line connecting distinct points on the boundary of the $D^2$, contradicting the assumption that $u$ is an immersion.

We therefore assume that the segment from $p$ to $x_{d}$ lies on the subset of the boundary of labeled by $\gamma_1$.  Since we chose $\gamma'_1$ to be $C^1$ close to $\gamma_1$, $x_d$ corresponds to an intersection point $x'_{d}$ between $\alpha_{n_d}$ and $\gamma_1'$.  The curve $\alpha_{n_d}$ separates $w$ into two triangles.  Let us denote by $v$ the triangle corners $p$, $x_d$ and $x'_d$.  In addition, we assume that $x'_d$ lies on the segment from $p$ to $x_{d+2}$ (If this is not true, we can simply reverse  the r\^oles of $x_d$ and $x_{d+2}$).

We think of $u$ as a map from a disc $D$ and consider the component $D'$ of the inverse image of $v$ under $u$ which contains the marked point $z_p$ corresponding to $p$. By our assumption, there are two points $z_1$ and $z_2$ on the boundary of $D'$, and segments from $z_p$ to $z_i$ along the boundary of the original disc $D$ which map to $\gamma_1 \cup \alpha_{n_d}$ and $\gamma'_1$ respectively.  The map $u$ only has critical points at $p$, $x_d$, and $x_{d+2}$, and at marked points that map to intersections between $\alpha_{n_i}$ and $\alpha_{n_d}$.  By $C^1$ closeness, none of these other possible marked points lie in $w$, so there can be no critical points on the remaining segment of $\partial D'$ which connects $z_1$ and $z_2$.  The map $u$ is therefore an immersion on this remaining segment, but this means that this segment must consist of a single point, i.e. that $z_1$ and $z_2$ are in fact the same point.  In particular, there are only three marked points, so we've proved the desired result.
\end{proof}

\begin{cor}
If $\beta$ is transverse to $\gamma_1$, and $\gamma'_1$ is as above, then the morphism $e_{1,1'}$ induces an isomorphism of chain complexes:
\[ m_{2}(e_{1,1'}, \_): CF^{*}(\beta, \gamma_1) \to   CF^{*}(\beta, \gamma'_1) ,\]
\end{cor}

Note that the proof of the above lemma did not require the assumption that $\Sigma$ is a cylinder.  We re-impose this condition in order to ensure that the curves are embedded.  The next result will allow us to conclude that $\gamma_1$ and $\gamma_2$ are also quasi-isomorphic.  For simplicity, we choose $\gamma'_1$ to intersect $\gamma_1$ in exactly two points which we denote $p$ and $q$.

\begin{lem}
Let $\gamma_1$, $\gamma'_1$ and $\gamma_2$ be closed embedded curves of equal holonomy, representing the generator of the first homology of the cylinder, such that $\gamma_1$ and $\gamma_2$ only bound embedded discs with disjoint interiors, and $\gamma_1$ and $\gamma'_1$ are $C^1$ close.  Up to possibly rescaling one of the terms by $\pm t^{\lambda}$ for some $\lambda$
\[ m_2(e_{2,1'}, e_{1,2}) = e_{1,1'} .\]
\end{lem}
\begin{proof}
We label the bigons $u_k$ with boundary $\gamma_1 \cup \gamma_2$ and their corners $p_i$ and $q_i$ as in the discussion following Corollary \ref{inductive_reduction}.  For the curves $\gamma'_1$ and $\gamma_2$ we also have corresponding bigons $u'_k$, and intersection points $p'_i$ and $q'_i$ which are $C^0$ close to their unprimed analogues. We must prove that
\[ \sum_{i,j} m_2( t^{ \sum_{k=0}^{2j-1} (-1)^{k} \omega(u'_k) } [q'_{j}],  t^{ \sum_{k=0}^{2i-1} (-1)^{k} \omega(u_k) } [p_i]) = \pm t^{\lambda} [p] .\]

In fact, the only term in the left hand side which does not vanish is the one corresponding to the corners of the bigon $u_k$ on whose boundary $p$ lies, see Figure \ref{unit_product}.  This proves the lemma.
\end{proof}

\begin{figure}
\begin{center}
\input{unit_product2.pstex_t} 
\end{center}
\caption{}

\label{unit_product}
\end{figure}

\begin{lem}
Let $\beta$ be an arbitrary unobstructed curve such that $\{\beta , \gamma_1, \gamma_2 \}$ is a transverse set.  Then
\begin{align*}
m_2(e_{1,2} ,\_ ) : CF^{*}(\beta, \gamma_1) & \to CF^{*}(\beta, \gamma_2) \\
m_2(\_, e_{1,2}) : CF^{*}(\gamma_2, \beta) & \to CF^{*}(\gamma_1, \beta)
\end{align*}
induce isomorphisms in homology. 
\end{lem}
\begin{proof}
We will only prove that $m_2(e_{1,2} ,\_ )$ is injective on homology, since the remaining parts of the lemma can be proved along the same line.
Let $\gamma'_1$ be $C^{1}$-close to $\gamma_1$.  Then, there is a bijection between $\gamma_1 \cap \beta$ and $\gamma'_{1} \cap \beta$, so that, $ CF^{*}(\beta, \gamma'_1)$ and $ CF^{*}(\beta, \gamma_1)$ are isomorphic as graded vector space. 

Let us now pass to homology, and observe that Lemma \ref{infty_relation} implies that $m_2$ descends to a product on homology which is associative.  The previous lemma implies that we can factor multiplication by $e_{1,1'}$ as multiplication by $e_{1,2}$ followed by multiplication by $e_{2,1'}$. We conclude that
\[ m_2(e_{1,2} ,\_ ) : CF^{*}(\beta, \gamma_1) \to CF^{*}(\beta, \gamma_2) \]
is injective on homology.
\end{proof} 
Note that this completes the proof for embedded curves in the cylinder.  We now prove the general case.
\begin{proof}[Proof of Proposition \ref{invariance} ]
Let $\beta$ be an arbitrary unobstructed curve such that $\{\beta , \alpha_1, \alpha_2 \}$ is a transverse set.  We need to show that
\[ m_2(e_{1,2}, \_ ) : CF^{*}(\beta, \alpha_1) \to CF^{*}(\beta, \alpha_2) \]
induces an isomorphism in homology.

For every $c \in  CF^{*}(\beta, \alpha_i)$ we note that the inclusion of $\alpha_1 \subset \Sigma$ induces a homomorphism
\[ \bZ \to \pi_{1}(\Sigma,c) \]
so we may consider the associated cover $\tilde{\Sigma}^{c}$ of $\Sigma$.  Consider a lift of $c$.  We obtain uniquely determined lifts $\tilde{\alpha}_{1}^{c}$ and $\tilde{\beta}^{c}$ since both $\alpha_1$ and $\beta$ pass through $c$.  Further, any of the intersection points $a_i$ between $\alpha_1$ and $\alpha_2$ determine the same lift $\tilde{\alpha}_2^{c}$  of our given curves (Other intersection points which don't appear in our morphism $e_{1,2}$ could determine other lifts, but that is of no concern to us). Further, $\tilde{\alpha}_{1}^{c}$, $\tilde{\alpha}_{2}^{c}$ are embedded circles in  $\tilde{\Sigma}^{c}$.  

We define the following equivalence relation on the intersection points $c$:
$c \cong c'$ if the unique basepoint preserving map from  $\tilde{\Sigma}^{c}$ to  $\tilde{\Sigma}^{c'}$ which is compatible with the projection to $\Sigma$ carries the curves $(\tilde{\alpha}_{1}^{c}, \tilde{\alpha}_{2}^{c} ,\tilde{\beta}^{c})$ to  $(\tilde{\alpha}_{1}^{c'}, \tilde{\alpha}_{2}^{c'} ,\tilde{\beta}^{c'})$.  Let 
\[ CF^{*}_{[c]}(\beta, \alpha_i) \]
denote the subgroup of $ CF^{*}(\beta, \alpha_i) $ which is generated by points in the equivalence class of $c$.  Since holomorphic discs lift to the universal cover,the decomposition
\[  CF^{*}(\beta, \alpha_i) = \bigoplus_{[c]}  CF^{*}_{[c]}(\beta, \alpha_i) .\]
is a decomposition into sub-complexes.  In fact, lifting holomorphic triangles also shows that $m_2 ( e_{1,2}, \_)$ respects this decomposition.  So it suffices to show that $m_2( e_{1,2}, \_)$ induces an isomorphism on homology for each sub-complex $CF^{*}_{[c]}(\beta, \alpha_i)$.  But 
\[ CF^{*}_{[c]}(\beta, \alpha_i)  \cong  CF^{*}(\tilde{\beta}^{c} , \tilde{\alpha}_{1}^{c}).\]
We are now in the situation of Floer homology for embedded curves in the cylinder, so the Proposition follows from the previous Lemma.
\end{proof}

\section{Cones and connect sums}
In this section, we consider the connection between twisted complexes and connect sums of curves.  For general symplectic manifolds, these results take the form of a long exact sequence for Floer homology \cite{seidel-les}.  For related consideration, see \cite{thomas} and \cite{TY}.
\begin{prop} \label{dehn-ext}
Let $\alpha$ and $\beta$ are unobstructed curves which intersect transversally and minimally, and assume further that $\beta$ is embedded.  If $\{ c_{i} \}_{i=1}^{n}$ and $\{b_{j} \}_{j=1}^{m}$ are the natural bases of $CF^1(\alpha,\beta)$ and $CF^{1}(\alpha, \beta[1])$ then the Dehn twist $\tau^{-1}_{\beta}(\alpha)$ is quasi-isomorphic to the twisted complex
\[
\xymatrix{  \alpha \ar[rrr]^-{\left([c_1], \cdots, [c_n], [b_1],\cdots , [b_m] \right)} & & & \beta^{\oplus n} \oplus \beta[1]^{\oplus m} . }
\] 

\end{prop}

This result follows from iteratively applying Lemma \ref{ext} which we will state and prove presently.  First, we need a definition.
\begin{defin}
If $\alpha$ and $\beta$ are immersed oriented curves which intersect at $c$, then $\alpha \#_{c} \beta$ is the immersed oriented curve which is obtained by resolving $c$ as in Figure \ref{ext-pict} and which satisfies
\[ \hol_A(\alpha \#_{c} \beta ) = \hol_{A}(\alpha) +  \hol_A(\beta).\] 
\end{defin}
\begin{rem}
The above equality simply says that the shaded triangles in Figure \ref{ext-pict} have equal area.  Note that the curve $\alpha \#_{c} \beta$ is only well defined up to isotopies which preserve area, and hence up to quasi-isomorphism in the Fukaya category by the results of the previous section.
\end{rem}

\begin{figure}
\begin{center}
\input{extension.pstex_t} 
\end{center}
\caption{}
\label{ext-pict}
\end{figure}
\begin{lem} \label{ext}
Let $\alpha$ and $\beta$ be unobstructed curves which intersect transversally and minimally, and let $c$ be an intersection point which has degree $1$ in $CF^{*}(\alpha, \beta)$.  Then $\alpha \#_{c} \beta$ is quasi-isomorphic in $Tw(\Fuk(\Sigma))$ to the twisted complex 
\[
\xymatrix{ Cone(c)= \alpha \ar[r]^-{[c]} & \beta. }
\] 
\end{lem}
\begin{proof}
Since the intersection point $a$ which appears in Figure \ref{ext-pict} has degree $0$ in  $CF^{*}(\alpha \#_{c} \beta, \alpha)$, we can consider it as an element of $CF^{0}(\alpha \#_{c} \beta, Cone(c))$.
\[ \xymatrix{ \alpha \#_{c} \beta \ar[d]^-{[a]} & \\
\alpha \ar[r]^-{[c]} & \beta. } \]
Proving that $[a]$ is a closed element is equivalent to showing:
\begin{align*}
m_1([a]) & = 0  
\\ m_2([c],[a]) & = 0 .\end{align*}
Since $\alpha$ and $\beta$ intersect minimally, there can be no bigon with boundary on $\alpha$ and $\alpha \#_{c} \beta $, so $m_1([a])$ necessarily vanishes. As for the second term, there are two triangles $u$ and $v$ that contribute to it, corresponding to the two shaded regions in Figure \ref{ext-pict}.  Our assumption that the holonomy of the connection $1$-form $A$ on $\alpha \#_c \beta$ is the sum of the holonomies on $\alpha$ and $\beta$ implies that these two triangles have the same area.  Neglecting the marked point, one checks that both discs are assigned a positive sign.  Since there are three marked points (one for each curve), the signed contributions cancel.  The same idea proves that $[b]$ defines a closed morphism in  $CF^{0}(Cone(c),\alpha \#_{c} \beta)$.
\[ \xymatrix{ \alpha \ar[r]^-{[c]} & \beta \ar[dl]^-{[b]} \\
\alpha \#_{c} \beta.  & } \]

Let $\gamma$ be an arbitrary curve which is transverse to $\alpha$ and $\beta$. We will prove that
\[ \hat{m}_2( [b], \hat{m}_2([a],\_)): HF^{*}( \gamma, \alpha \#_c \beta) \to  HF^{*}( \gamma, \alpha \#_c \beta) \]
is an isomorphism, which establishes that $m_2([a], \_)$ is injective on homology.  One should check its surjectivity, as well as the isomorphism of $m_2(\_,[a])$, but the proofs are entirely analogous, so we shall omit them.  Note that the following lemma clearly implies the desired result.
\end{proof}
\begin{lem}
If $\alpha \#_{c} \beta$ is $C^{1}$ close to $\alpha \cup \beta$ away from a small neighbourhood of $c$ then
\[\hat{m}_2( [b], \hat{m}_2([a],\_)):  CF^{*}( \gamma, \alpha \#_c \beta) \to  CF^{*}( \gamma, \alpha \#_c \beta) \]
is an isomorphism of $\Lambda$ modules.
\end{lem}
\begin{proof}
The condition of $C^{1}$ closeness implies that, to an intersection point $x$ between $\gamma$ and  $\alpha \#_{c} \beta$, there corresponds a unique closest intersection point $x'$ between $\gamma$ and $\alpha$ or $\beta$.  Assume that $x'$ lies on $\beta$ (the case where it lies on $\alpha$ is similar).  Note that the segment between $x$ and $x'$ separates the triangle $u$ into a square $u_1$ with corners $x$, $a$, $c$ and $x'$, and a triangle $u_2$ with corners $b$, $x'$ and $x$.  In the simplest case, these are the only polygons that contributes to the composition, and we easily compute
\begin{align*}
m_3([c],[a],[x]) & = \pm  t^{\omega(u_1)} [x'] \\
m_2([b],[x']) & = \pm t^{\omega(u_2)} [x]
\end{align*}
which is the desired result.  However, there may be other polygons.  To understand why these other polygons do not affect out result, we observe that $C^1$ closeness implies that any other polygons must have area larger than the area of $u$.  Indeed, by passing to the universal cover, it is easy to see that any such polygons must either wrap around $\alpha$ or $\beta$ more than once, in which case they clearly have area larger than the area of $u$, or they must have corners which include an intersection point which is close to $\alpha$, and one which is close to $\beta$, which means that such a polygon has area approximately equal to the area of a corresponding polygon with one corner at $\alpha$, the other at $\beta$, and the last corner being $c$.  It is clear that we can choose $u$ to have area arbitrarily smaller than the area of all such polygons.

This implies that the lowest order term in the determinant of the matrix of the composition $\hat{m}_2( [b], \hat{m}_2([a],\_))$ is given exactly by
\[ \pm t^{\omega(u) |\gamma \cap \alpha \#_{c} \beta|} \]
which is clearly non-zero.
\end{proof}

\section{Trivial objects in the derived category} \label{triv-obj}
We would like to define homomorphisms from $K(\Fuk(\Sigma))$ to the groups mentioned in the Introduction.  In order to do this, we choose a vector field $\nu$ which vanishes at only one point (see Appendix \ref{wind_holo}). Consider the map
\[ \bZ^{|\Ob(\Fuk(\Sigma))|} \to H_{1}(\Sigma, \bZ) \oplus \bR \oplus \bZ / \chi(\Sigma) \bZ \]
defined by
\[ \Gamma = (\gamma_k) \mapsto ([\Gamma], \hol_{A}(\Gamma), \wind_{\nu}(\Gamma)) \equiv \sum_{k} ( [\gamma_k] , \hol_{A}(\gamma_k) , \wind_{\nu} (\gamma_k) ) ,\]
where $\Gamma$ is an arbitrary twisted complex over the Fukaya category of $\Sigma$.  Recall that this means that we have a finite collection of unobstructed curves $ \{ \gamma_{k} \}$ and morphisms of degree $1$, $\Delta=(\delta_{k,l})$ for $k<l$ such that $\hat{m}_1(\Delta)=0$.  In this section, we prove

\begin{prop} \label{zero}
If $\Gamma$ is quasi-isomorphic to the zero object then
\[  ([\Gamma], \hol_{A}(\Gamma), \wind_{\nu}(\Gamma)) =0 .\]
\end{prop}

We begin with some preliminary remarks.  For a fixed $k$, consider the set of intersections $\cup_{l} (\gamma_k \cap \gamma_l)$.  This is a finite set of points on the curve $\gamma_k$ which subdivides it into segments.  Let  $\cI^{\gamma_k}$ denote the set of such segments.  Let us choose $C^1$ small Hamiltonian perturbations $\gamma'_k$ of the curves $\gamma_k$ satisfying the following conditions:
\begin{itemize}
\item The curves $\gamma_k$ and $\gamma'_k$ intersect twice in each segment $I \in \cI^{\gamma_k}$.  Further, if we orient $I$ according to the orientation of $\gamma_k$, then the intersection point which has degree $0$ in $CF^*(\gamma_k, \gamma'_k)$ occurs before the intersection point which has degree $1$.
\item There is a constant $\lambda$ independent of $k$ such that all the bigons bounded by $\gamma_k$ and $\gamma'_k$ have area $\lambda$.
\item Consider a local model for an intersection point $c_{k,l}$ between $\gamma_k$ and $\gamma_l$.  There is a constant $\lambda'$ independent of $k$, $l$, and the intersection point $c_{k,l}$  such that all the triangles shaded in Figure \ref{perturbation} have the same area.
\end{itemize}

\begin{figure}
\begin{center}
\input{perturbation.pstex_t}
\end{center}
\caption{}

\label{perturbation}
\end{figure}

We decompose the differentials $\delta_{k,l}$ in term of the generators corresponding to the intersections between $\gamma_k$ and $\gamma_l$
\[ \delta_{k,l} = \sum_{c_{k,l} \in \gamma_{k} \cap \gamma_{l}} P_{c_{k,l}} [c_{k,l}] .\]
Of course, by definition, $P_{c_{k,l}}$ vanishes whenever $c_{k,l}$ does not have degree $1$ in $CF^*(\gamma_k, \gamma_l)$.  Our assumption of $C^{1}$ closeness implies that there is a unique intersection point $c'_{k,l}$ between $\gamma'_k$ and $\gamma'_l$ corresponding to $c_{k,l}$.  We define
\[ \delta'_{k,l} = \sum_{c'_{k,l} \in \gamma'_{k} \cap \gamma'_{l}} P_{c_{k,l}} [c'_{k,l}] .\]
Let $\Gamma'$ be the collection of curves $\gamma'_k$, and $\Delta'$ be the collection of morphisms $\delta'_{k,l}$.
\begin{lem}
$(\Gamma', \Delta')$ is a twisted complex.
\end{lem}
\begin{proof}
Each $c_{k,l}$ that contributes to $\delta_{k,l}$ has degree $1$, so Figure \ref{perturbation} is an accurate  model for the input corners of each polygon that contributes to $\hat{m}_{1}(\Delta)$.  In addition, since the output of $\hat{m}_1(\Delta)$ has degree $0$, this is also an accurate model for the output if we let the horizontal curve be the target, and the vertical curve be the source. Note that in the local model near every intersection point, each polygon $u$ which contributes to $\hat{m}_1(\Delta)$ lies in the second or fourth quadrant.

Such a polygon can be uniquely perturbed to a polygon $u'$ which contributes to $\hat{m}_{1}(\Delta')$.  We simply go along the curves $\gamma'_k$ instead of $\gamma_k$.  Our assumption that all the triangles in Figure \ref{perturbation} have the same area implies that $u$ and $u'$ have the same area. This can be proved by considering a neighbourhood of each corner $c_{k,l}$.  For example, if $u$ locally lies in the second quadrant, then $u'$ contains the bigon that lies on the negative $x$-axis, but misses the union of two triangles and a square.  The bigon has area $\lambda$, but this is exactly the sum of the area of the two triangles and the square.  It is clear that the respective signed contributions of $u$ and $u'$ to $\hat{m}_1(\Delta)$ and $\hat{m}_{1}(\Delta')$ are equal, so the vanishing of the first implies the vanishing of the second.
\end{proof}

Let $e_k$ be the quasi-isomorphism of $CF(\gamma_k, \gamma'_k)$ constructed in Lemma \ref{invariance}.  Since all the bigons bounded by $\gamma_k$ and $\gamma'_k$ have the same area, we can set
\[ e_k = \sum_{p \in |\gamma_k \cap \gamma'_k| \ni \deg(p)=0} [p] .\]
We define
\[ E = \oplus_{k} e_k \in \bigoplus_{k} CF^{*}(\gamma_k, \gamma'_k) \subset CF^{*}(\Gamma, \Gamma').  \]

\begin{lem}
$E$ is a closed morphisms of twisted complexes.  
\end{lem}
\begin{proof}
Since the chosen perturbation is $C^1$-small the argument of Lemma \ref{strict_units} implies that the map
\[ m_2(e_k, \_): CF^{*}(\gamma_l, \gamma_k) \to CF^{*}(\gamma_l, \gamma'_k) \]
is an isomorphism of vector spaces, and hence has an inverse.  Further, since $C^1$ closeness of $\gamma'_i$ and $\gamma_i$ relative to $\gamma_j$ for $i \neq j$ implies $C^1$ closeness relative $\gamma'_j$,  Lemma \ref{strict_units} implies that the higher products $m_{l+1}(x'_{l}, \cdots, x'_k, e_k, x_{k-1}, \cdots, x_1 )$ vanish whenever $l>1$ and all the morphisms other than $e_k$ are either amongst curves $\gamma_j$ for $j \leq k$ or $\gamma'_j$ for $j \geq k$. In particular, proving that $\hat{m}_1(E)=0$ is equivalent to
\begin{align*}
m_1(e_k) &= 0
\\ m_2(e_l, \delta_{k,l}) + m_2(\delta'_{k,l}, e_k) &= 0.
\end{align*}
The first condition only involves the curves $\gamma_k$ and $\gamma'_k$, and hence holds by construction.  The second condition follows from observing that for each $p_l$, there exists at most one $c_{k,l}$ that contributes to the first term of the above equation.  This intersection point $c_{k,l}$ must be one of the endpoints of the segment of $\gamma_l$ on which $p_l$ lies.  The consistent choice of area of bigons implies that
\[ m_2(p_l, c_{k,l}) + m_2(c'_{k,l}, p_k) = 0 \]
whenever $c_{k,l}$ is a morphism that contributes to $\delta_{k,l}$ and $p_k$ is the degree $0$ intersection between $\gamma_k$ and $\gamma'_k$ which is adjacent to $c_{k,l}$.  The argument of Lemma \ref{strict_units} implies that no other degree $0$ intersection between $\gamma_k$ and $\gamma'_k$ pairs non-trivially with $c'_{k,l}$.  So the computation reduces to the local model of Figure \ref{perturbation}.  We leave it to the reader to check the local computation, and simply comment that one triangle of area $\lambda'$ contributes to each of the its terms, and that the only thing to check is that they have opposite signed contributions.
\end{proof}
\begin{rem}
Using these ideas, one can easily prove that $E$ is a quasi-isomorphism. We do not need this fact.  Also, one can approach these constructions from a more algebraic point of view.  Using the Lemma \ref{invariance}, one can formally construct a new twisted complex $\Gamma'$ that is quasi-isomorphic to $\Gamma$.  
\end{rem}

Let $\cK$ be the collection of increasing indices among the integers that label the curves $\gamma_k$.  In the notation of Appendix \ref{twisted-and-K}, the fact that $(\Gamma, \Delta)$ is a twisted complex corresponds to the vanishing of 
\[ \sum_{K \in \cK} m_{|K|}(\Delta_{K}) .\]
We write $\cC$ for the collection of all sequences $( c_{k_{n},k_{n+1}} )$ with $k_{n} < k_{n+1}$.  Every element of $\cC$ is labeled by some increasing sequence $K \in \cK$.  The vanishing of $\hat{m}_1(\Delta)$ is equivalent to
\[ \sum_{C_{K} \in \cC} P_{C_K} m_{|K|}(C_K) = 0 ,\]
where 
\[ P_{C_K} = \prod_{c_{k,l} \in C_K} P_{c_{k,l}} .\]
We will sometimes write $C$ for $C_K$ when the notation becomes sufficiently cumbersome.  We therefore have
\begin{lem}
The sets of polygons 
\[ \cM(g, C_{K}) = \cM(g, c_{k_{|K|-1},k_{|K|}}, \ldots, c_{k_0, k_1} ) \]
with $C_K \in \cC$ are exactly those that contribute to $\hat{m}_{1}(\Delta)$.
\end{lem} \noproof

The assumption that $\Gamma$ is quasi-isomorphic to the $0$-object implies that there exists $H \in CF^{*}(\Gamma, \Gamma')$ such that $\hat{m}_{1}(H) = E$. As above, the sets of polygons that contribute to $\hat{m}_1(H)$ are
\[ \cM(x  , C'_{K_1} ,h ,C_{K_2}) .\]
If $x = p_k$, then the first element of $K_1$ and the last element of $K_2$ agree.  In particular, since $K_i$ are both increasing sequences, $h$ must lie in $CF^{1}(\gamma_l, \gamma'_k)$ where $k \leq l$ in order to contribute to 
\[ \langle \hat{m}_1(H) , [p_k] \rangle = 1 ,\]
where the left hand side is the coefficient of $[p_k]$ in $\hat{m}_1(H)$.  Further, since the differential vanishes identically in $CF^{1}(\gamma_l, \gamma'_l)$, we may assume that we have strict inequality.  Let $F$ be the sum of the terms in $H$ which lie in $CF^1(\gamma_l, \gamma'_k)$ for $k < l$.  We write
\[ F= \sum_{f \in \gamma_l \cap \gamma'_k | \deg(f) =1} P_{f} [f] .\]
  By construction, we still have
\begin{equation}
\langle \hat{m}_1(F), [p_k] \rangle = 1 \end{equation} 
for every degree $0$ intersection $p_k$ between $\gamma_k$ and $\gamma'_k$.

Again, there is a canonical bijection between $\gamma_l \cap \gamma'_k$ and $\gamma_k \cap \gamma_l$.  Let $\bar{f}$ be the intersection point between $\gamma_k$ and $\gamma_l$ corresponding to $f$. $[\bar{f}]$ is a generator of degree $0$ in $CF^{*}(\gamma_k, \gamma_l)$.  The next lemma, though seemingly innocuous, is the key to our arguments in this section.

\begin{figure}
\begin{center}
\input{perturbed_square.pstex_t}
\end{center}
\caption{}

\label{perturbed_square}
\end{figure}

\begin{lem} \label{polygon_correspondance}
If $k \in K$ and the sequences $K_1$ and $K_2$ consist respectively of the elements which are smaller and larger than $k$, then every polygons $u$ in
\[ \cM(\bar{f}, C_K) \]
corresponds to a finite set of polygons $\{v_r \}$ in
\[ \cM(p_k , C'_{K_1} , f , C_{K_2}) .\]
Further, all the polygons corresponding to $u$ have the same area which differs from that of $u$ by $\lambda-\lambda'$.
\end{lem}
We will denote this finite number of polygons by $|n_u(p_k)|$.  The notation is justified by
\begin{defin}
If $u$ is a polygon with a boundary edge on $\gamma$, and $p$ is a point on $\gamma$, then $n_u(p)$ is the signed multiplicity of $p$ as a point on the edge of $\partial u$ which lies on $\gamma$. \end{defin}
Thinking of $u$ as an orientation preserving map from the disc, there is a segment on $\partial D$ labeled by $\gamma$.  The non-negative number $|n_{u}(p)|$ is the number of points on this edge of $\partial D$ which are mapped to $p$.  The sign is positive if and only if the orientation of $\gamma_k$ agrees with that of boundary of $D$ at $p$.
\begin{proof}[Proof of Lemma \ref{polygon_correspondance}]
Let us assume for simplicity that $n_u(p_k)=1$.  We can perturb every $u$ in $\cM(\bar{f}, C_k)$ to a polygon $v$ whose boundary acquires a new convex corner between $\gamma_k$ and $\gamma'_k$ at $p_k$.  Because $p_k$ is our output, as we are moving along the boundary of the polygon, we are in fact moving from $\gamma'_k$ to $\gamma_k$ at $p_k$.  The boundary of the polygon then continues along the curves $\gamma_l$ for $k < l$ going through the corners which appear in $C_{K_2}$.  Eventually, we arrive in a neighbourhood of $f$.  While $u$ has a corner at $\bar{f}$ and picks up a boundary segment along $\gamma_{k_1}$ where $k_1$ is the first element of $K$, $v$ has a corner at the nearby intersection point $f$, and its next boundary segment is along $\gamma'_{k_1}$.  The remaining segments of the boundary of $v$ lie on the curves $\gamma'_l$ for $l<k$ with corners at $C'_{K_1}$.

In order to prove that the difference between the area of $u$ and $v$ is $\lambda - \lambda'$, one computes this difference locally in a neighbourhood of each corner as in the local model of Figure \ref{perturbation}.  Instead of writing the argument, we illustrate the case with four corner in Figure \ref{perturbed_square}.  The dark (red) regions  lie in $u$ but not in $v$, and vice versa for the light (blue) regions.  The areas of these regions cancel pairwise except for the top right corner.  By construction, it has area $\lambda - \lambda'$.  One can perform the same computation if the multiplicity in negative, as in Figure \ref{perturbed_square_opposite}.

If $|n_u(p_k)|>1$, then the boundary of $u$ wraps around $\gamma_k$ more than once.  Say for specificity that $n_u(p_k) = 2$.  Then we can choose our perturbed polygons to either switch from $\gamma'_k$ to $\gamma_k$ at $p_k$ then wrap around $\gamma_k$ once, or to wrap around $\gamma'_k$ before switching at $p_k$ to $\gamma_k$.  The difference between the areas of these two polygons is exactly the area between $\gamma_k$ and $\gamma'_k$, which vanishes by assumption.  
\end{proof}

\begin{figure}
\begin{center}
\input{perturbed_square_opposite.pstex_t}
\end{center}
\caption{}
\label{perturbed_square_opposite}
\end{figure}
We will use the sets of polygons of the above Lemma to construct a $2$-chain $U$ whose boundary is $\sum_{k} [\gamma_k]$.  First, we write,
\[ P_{f} = \sum_{j} b^{f}_{j} t^{\beta^{f}_{j}} .\]
Similarly, we can decompose
\[ P_{c_{k,l}} = \sum_{i_{k,l} \in \bZ} a^{c_{k,l}}_{i_{k,l}} t^{\alpha^{c_{k,l}}_{k,l}} .\]
We can therefore compute
\begin{align*}
 P_{f} \langle  \hat{m}_1(\Delta), \bar{f} \rangle & = \sum_{C \in \cC, u \in \cM(\bar{f},C)} P_{C} P_{f} (-1)^{s(u)} t^{\omega(u)} \\
& = \sum_{C \in \cC, i, j, u \in \cM(\bar{f},C)} (-1)^{s(u)} a^{C}_{i} b^{f}_{j} t^{\omega(u) + \alpha^{C}_i + \beta^{f}_{j}},
\end{align*}
where $P_{C}$ is the product of the Laurent series $P_{c_{k,l}}$ and is decomposed as
\[ \sum_{i} a^{C}_{i} t^{\alpha^{C}_{i}} .\]
Note that $a^{C}_{i}$ is a products of monomials  $a^{c_{k,l}}_{i_{k,l}}$ for all $c_{k,l}$ that appear in $C$, and $\alpha^{C}_{i}$ is the sum of the corresponding exponents  $\alpha^{c_{k,l}}_{i_{k,l}}$.  The vanishing of $\hat{m}_1(\Delta)$ is equivalent to the vanishing of
\begin{equation} \label{closed_diff} \sum_{C , u \in  \cM(\bar{f},C) | \omega(u) + \alpha^{C}_{i} = \tau} (-1)^{s(u)} a^{C}_{i} \end{equation}
for fixed $\bar{f}$ and for any $\tau$.  This will be particularly useful if we let $\tau = b_{j}^{f}$ for some $j$.

There are only finitely many polynomials that contribute to the $t^{0}$ component of $P_{f} \langle \hat{m}_1(\Delta), \bar{f} \rangle$.  We introduce sets of polygons:
\begin{align*} \cU_{f,C,i,j} & = \{ u \in \cM(\bar{f},C) | \omega(u) + \alpha^{C}_i + \beta^{f}_{j} = -(\lambda - \lambda') \} \\
\cU & = \bigcup_{f,C,i,j} \cU_{f,C,i,j} .
\end{align*}
Note that every $u \in \cU$ determines $f$ and $C$ which essentially correspond to its corners, but not $i$ and $j$.  We define integers
\[ M_{i,j}(u) = (-1)^{s(u)} a^{C}_{i} b^{f}_{j} \]
whenever $u \in \cU_{f,C,i,j}$ and their sum
\[ M(u) = \sum_{i,j} M_{i,j}(u).\]
This allows us to consider the finite $2$-chain
\[ U = \sum_{u \in \cU} M(u) [u] \]
where $[u]$ carries the restriction of the orientation of $\Sigma$.

\begin{lem} \label{bounding_chain}
$U$ is a $2$-chain with boundary $\sum_k [\gamma_k]$.
\end{lem}
\begin{proof} Pick a generic point $p$ on $\gamma_k$.  We must show that the (weighted) number of discs which pass through $p$ is precisely one.  In other words, we must show that:
\begin{equation}\label{degree1} \sum_{u \in \cU} M(u) n_{u}(p) = 1 . \end{equation} 
Note that $p$ lies on a segment $I \in \cI^{\gamma_k}$ and hence lies on the same segment a unique generator $p_k$ of $CF^{0}(\gamma_k, \gamma'_k)$.  Since every disc passing through $p$ also passes through $p_k$, we may assume that $p=p_k$ for some degree $0$ intersection $p_k$ between $\gamma_k$ and $\gamma'_k$.

It will suffice to show that the left hand-side of Equation \eqref{degree1} is the coefficient of $p_k$ in the Novikov series expansion of the left hand-side of
\[ \hat{m}_1(H) = E .\]

Except for the signs, this is the content of Lemma \ref{polygon_correspondance}. Indeed, the conditions \[ \omega(u) + \alpha^{C}_i + \beta^{f}_{j} = -(\lambda - \lambda') \]
implies the the discs $v_r$ satisfy
 \[ \omega(v_r) + \alpha^{C}_i + \beta^{f}_{j} = 0 \]
 The proof would therefore be complete if we check that the sign that each $v_r$ in the statement of Lemma  \ref{polygon_correspondance} contributes to $\hat{m}_1(H)$ is $(-1)^{s(u)+s'(u)}$ where $s(u)$ is the sign that $u$ contributes to $\hat{m}_1(\Delta)$ and which already appears in Equation \eqref{degree1}, and $s'(u)$ is a correction which vanishes $\mod 2$ if and only if the given orientation of $\gamma_k$ and its induced orientation as the boundary of $u$ agree at $p_k$.

Since $p_k$ has degree $0$ in $CF^*(\gamma_k, \gamma'_k)$, it does not affect the sign of the contribution of $v_r$.  Also, every intersection point which appears in $C$ will contribute the same sign to $\hat{m}_1(H)$ and $\hat{m}_1(\Delta)$.  Therefore, $s'(u)$ is just the difference between the contribution of $\bar{f}$ to the sign of $u$ and of $f$ to the sign of $v_r$.  But $\bar{f}$ has vanishing degree, so $s'(u)$ is just the contribution of $f$ to the sign of $v_r$.  This contribution vanishes if and only if the orientation of $\gamma'_{k_1}$ at $f$ agrees with its induced orientation as the boundary of $v_r$.  Since every intersection point which appears in $C$ has degree $1$, agreement of orientations at $f$ implies agreement of orientations at $p_k$, which was the desired result (see Figure \ref{perturbed_square} for one of the two possible cases when considering squares).
\end{proof}

This proves the first part of Proposition \ref{zero}.  We now prove the second part.

\begin{lem}
The signed area of the 2-chain $U$ vanishes.
\end{lem}
\begin{proof}
We must prove that
\begin{equation} \label{area_series} \sum_{u \in \cU} \sum_{i,j} M_{i,j}(u) \omega(u) = 0.\end{equation}

We replace $\omega(u)$ by $ \alpha^{C}_i + \beta^{f}_{j} - (\lambda - \lambda')$, and observe that, following Equation \eqref{closed_diff}, the vanishing of the $[\bar{f}]$ component of $\hat{m}_1(\Delta)$ implies that for fixed $j$
\[  \sum_{u \in \cU_{f},i} M_{i,j}(u) = 0 ,\]
so the term containing $\beta^{f}_{j}-(\lambda - \lambda')$ vanishes. We now recall that $\alpha^{C}_i$ is itself the sum of contributions $\alpha^{c_{k,l}}_{i_{k,l}}$.  Our goal is to show that the total coefficient of each $\alpha^{c_{k,l}}_{i_{k,l}}$ vanishes.  

For simplicity, let us assume that we have a monomial
\[ P_{c_{k,l}} = a^{c_{k,l}} t^{\alpha^{c_{k,l}}} .\]
In this case, $\alpha^{c_{k,l}}$ appears in every exponent $\alpha^{C}_{i}$ whenever $c \in C$, so we will want to show
\begin{equation} \label{vanishing_at_intersection} \sum_{ u \in \cU | c_{k,l} \in\partial^{0} u} M(u) =0 ,\end{equation} 
where $\partial^{0}u$ are the corners of $u$.  Topologically, this expression corresponds to the signed multiplicity of the chain $U$ in the four quadrants surrounding $c_{k,l}$.  By the previous Lemma, the total sum therefore vanishes.

The general case where $P_{c_{k,l}}$ is not a monomial is addressed similarly.  Indeed, assume that $u \in \cU$ contains $c_{k,l}$ as a corner, and define $C' = C - \{c_{k,l} \}$.  Just as we attached a Laurent series to $C$, we define
\[ P_{C'} = \sum_{i'} a^{C'}_{i'} t^{\alpha^{C'}_{i'}} \]
as the product of the Laurent series associated to each corner in $C'$.  Note that the coefficient of $ a^{c_{k,l}}_{i_{k,l}} \alpha^{c_{k,l}}_{i_{k,l}}$ in Equation \eqref{area_series} is
\[ \sum_{u \in \cU, i', j |c_{k,l} \in \partial^0(u) } (-1)^{s(u)} b^{f}_{j} a^{C'}_{i'} \]
where the corner locus of $u$ consists of $c_{k,l}$, $\bar{f}$, and the elements of the set $C'$, while the coefficients $i'$ and $j$ are required to satisfy
\[ \beta^{f}_{j} + \alpha^{C'}_{i'} + \omega(u) + (\lambda - \lambda')= - \alpha^{c_{k,l}}_{i_{k,l}} .\]
This series is therefore exactly  coefficient of  $t^{-\alpha^{c_{k,l}}_{i_{k,l}}}$ in 
\[ \sum_{u \in \cU, c_{k,l} \in \partial^{0}u} (-1)^{s(u)} t^{\omega(u) + (\lambda - \lambda')} P_f P_{C'} .\]
Our goal is to show that this series vanishes.  Since the Novikov ring has no zero divisors, there is not harm in multiplying the above with $P_{c_{k,l}}$.  But the resulting polynomial is 
\[ \sum_{u \in \cU, c_{k,l} \in \partial^{0}u} (-1)^{s(u)} t^{\omega(u)+ (\lambda - \lambda')} P_f P_C .\]
Using the proof of the previous Lemma, we see that this is exactly the difference between the coefficient of $p^0_k$ and $p^1_{k}$ in $\hat{m}_1(H)$ and hence vanishes.
\end{proof}

\begin{cor} \[\hol_A(\Gamma) = 0\]

\end{cor}
\begin{proof}
For an immersed polygon $u$, we have
 \[\hol_A(\partial u) =  \omega(u).\]
The result therefore follows from Lemma \ref{bounding_chain} and the previous lemma.  For details, the reader should go through the proof of the next lemma, replacing the winding number with the holonomy of the connection.
\end{proof}

We now complete the proof of Proposition \ref{zero}.
\begin{lem}
\[ \wind_{\nu}(\Gamma) = 0 \]
\end{lem}
\begin{proof}
The strategy is to relate the winding number around the curves $\gamma_i$ to the winding number around the boundaries of polygons $u \in \cU$. First, we observe that
\[
\sum_{u \in \cU} M(u) \wind_{\nu}(\partial u)  = 0
\] 
since $ \wind_{\nu}(\partial u) = -1$ and the above sum of coefficients therefore vanishes by Equation \eqref{closed_diff}.  We can summarize this by writing
\[ \wind_{\nu} (\partial U) = 0 .\]

On the other hand, we can compute the winding number by using the decomposition of $\partial u$ into piecewise smooth curves which lie on the curves $\gamma_i$ as
\[ \wind_{\nu}(\partial u) = \sum_{c_{k,l} \in C} \wind_{\nu}(I(\gamma_{k},\partial u
)) - \theta(c_{k,l}),\]
where $ I(\gamma_{k},\partial u) $ is the segment between two successive corners of $u$ along $\gamma_{k}$ oriented as a boundary of $u$, and $\theta(c_{k,l})$ is the angle between the two boundary segments of $u$ that share $c_{k,l}$ as a vertex.  Because we're not taking the orientations of $\gamma_{k}$ and $\gamma_{l}$ at $c_{k,l}$ into account, this angle $\theta(c_{k,l})$ is the same for all polygons that have a corner at $c_{k,l}$ and contribute to $\hat{m}_1 (\Delta)$.

Using this formula, we compute
\begin{equation*}
\wind_{\nu}(\partial U)  =  \sum_{c_{k,l}} -\theta(c_{k,l}) \sum_{u \in \cU, c_{k,l} \in \partial^{0} u} M(u) + \sum_{\gamma_k} \sum_{u \in \cU}  M(u) \wind_{\nu}(I(\gamma_k,\partial u)).
\end{equation*}

The first series vanishes as in the previous Lemma.  As for the second series, we note that $I(\gamma_k, \partial u)$ may not have the same orientation as $\gamma_k$.  But the sign difference is $s'(u)$ which we introduced in the proof of Lemma \ref{bounding_chain}.  After decomposing  $I(\gamma_k, \partial u)$ into minimal segments $I$ we see that the signed contribution of each such minimal segment is weighted by $n_{u}(p^I_{k})$,  where $p^I_k \in I$ is the degree $0$ intersection between $\gamma_k$ and $\gamma'_k$ that occurs in $I$.  So, 

\[ \wind_{\nu}(\partial U) =  \sum_{\gamma_k \in |\Gamma|} \sum_{I \in \cI^{\gamma_k}} \sum_{u \in \cU} M(u) n_{u}(p^I_k) \wind_{\nu}(I). \]
Note that in our proof that $\partial U = \sum [\gamma]$, we proved that for a fixed $I \in \cI^{\gamma_k}$,
\[   \sum_{u \in \cU} M(u) n_{u}(p^I_k) = 1. \]
In particular, we may now conclude
\begin{align*}
\wind(U) & =  \sum_{\gamma_k \in |\Gamma|} \sum_{I \in \cI^{\gamma_k}} \wind_{\nu}(I) \\
& = \sum_{\gamma_k \in |\Gamma|} \wind_{\nu}(\gamma_k),
\end{align*}
thereby completing the proof that this last sum vanishes.
\end{proof}

\section{Relations from the Mapping Class Group} \label{relations}
We will prove the following
\begin{prop} \label{homologous}
If $\alpha_1$ and $\alpha_2$ are homologous embedded curves such that 
\[ (\hol_A(\alpha_1), \wind_{\nu}(\alpha_1)) =  (\hol_A(\alpha_2), \wind_{\nu}(\alpha_2)) \]
then $\alpha_1$ and $\alpha_2$ represents the same class in $K(\Fuk(\Sigma))$.
\end{prop}
\begin{rem}
The reader may wonder what value to give a pictorial proof along the lines of Figure \ref{pair_pants_torus_proof}.  Such techniques are standard in the proofs of theorems about isotopy classes of curves, whereas we need to keep track of area.  But we already proved in Lemma \ref{invariance} that any two isotopic curves for which the holonomies of the connection $A$ are equal must be quasi-isomorphic, so our proofs are valid.
\end{rem}

We will need several lemmata in order to prove this result.  We begin by obtaining the first non-trivial relation in $K(\Fuk(\Sigma))$.
\begin{lem}
If $\gamma_1$ and $\gamma_2$ are unobstructed simple closed curves which bound symplectomorphic submanifolds of $\Sigma$, then they represent the same class in $K(\Fuk(\Sigma))$.
\end{lem}
\begin{proof}
There exists an element $\phi$ of the Hamiltonian mapping class group which maps $\gamma_1$ to $\gamma_2$.  Since we can factor $\phi$ as a product of symplectic Dehn twists, it suffices to show that if $\alpha$ is any simple closed curve on $\Sigma$, then $\tau_{\alpha}(\gamma_1)$ represents the same element of $\gamma_1$ in $K$-theory.  Since $\gamma_1$ is separating, the intersection number $\gamma_1 \cdot \alpha$ vanishes for any curves $\alpha$, so the result follows from Proposition \ref{dehn-ext}.
\end{proof}

This lemma easily generalizes as follows
\begin{lem}
If $S_1$ and $S_2$ are unions of simple closed curves which bound symplectomorphic submanifolds of $\Sigma$, then $S_1$ and $S_2$ represent the same class in $K(\Fuk(\Sigma))$. \noproof
\end{lem}
Note that this lemma implies that given any real number $\Omega$ less than the area of $\Sigma$, there is a well defined class in $K$-theory which represents the difference between isotopic curves which form the boundary of a embedded cylinder of area $\Omega$.  

\begin{lem} \label{cylinders}
The subgroup of $K$-theory generated by boundaries of embedded cylinders is isomorphic to $\bR$.
\end{lem}
\begin{proof}
By restricting the holonomy map to this subgroup, we have a surjection to $\bR$.  To prove that this map is injective, we consider an element of $K$-theory mapping to $0$ under the holonomy map.  Without loss of generality, all the terms are boundaries of cylinders that have the property that the absolute value of their area is less than $\frac{\Omega_{\Sigma}}{2}$.  Since the total holonomy vanishes, there must be at least two terms, one of which has positive, and the other negative holonomy.  By the previous Lemma we can represent the sum of these two terms as the boundary of one cylinder, which still has area at most $\frac{\Omega_{\Sigma}}{2}$.  Proceeding inductively, we end up with the boundary of a cylinder of vanishing area.  The corresponding element of $K$-theory vanishes by Lemma \ref{invariance}.  
\end{proof}

We denote the class of the boundary of a cylinder of area $\Omega$ by $\Omega \rho$.  Also, the class of a curve bounding a torus of zero area is a well defined element in $K(\Fuk(\Sigma))$.  We will denote this element by $T$.

\begin{lem}
If $\alpha_1$, $\alpha_2$ and $\alpha_3$ bound a pair of pants of area $\Omega$, then, in $K$-theory
\[ [\alpha_1] +  [\alpha_2] + [\alpha_3] = T + \Omega \rho\]
\end{lem}
\begin{figure}
\begin{center}
\input{pair_pants_torus.pstex_t}
\end{center}
\caption{}
\label{pair_pants_torus}
\end{figure}
\begin{proof}

Assume none of the curves is separating.  If the genus of $\Sigma$ is more than $2$, then we can find a curve $\alpha_4$ as in Figure \ref{pair_pants_torus}.  The intersections between $\alpha_4$ and $\alpha_2$ are essential since otherwise $\alpha_3$ bounds a torus.

Since we are not in the genus 2 case, there must be an additional curve $\gamma$ as in Figure \ref{pair_pants_torus_proof}. The left hand side corresponds to starting with $\alpha_1[1]$, then successively ``adding'' $\gamma$ and $\alpha_4$, whereas the right hand side starts with $\gamma$ and ``adds'' $\alpha_2$ and $\alpha_3$ to it. Since the final curves of the left and right hand columns of Figure \ref{pair_pants_torus_proof} are isotopic, we conclude that,
\[ [\alpha_2] + [\alpha_3] + [\gamma] = [- \alpha_1] + [\alpha_4] + [\gamma] \]
as long as $\alpha_4$ bounds a torus of appropriate area.  But additivity of the holonomy map implies that the appropriate area is the area of the pair of pants bounded by $\alpha_1$, $\alpha_2$ and $\alpha_3$.  This establishes the desired result.
\begin{figure}
\begin{center}
\input{pair_pants_torus_proof.pstex_t}
\end{center}
\caption{}
\label{pair_pants_torus_proof}
\end{figure}
In the genus 2 case, if the curves are non-separating, then we are in the situation of Figure \ref{genus_2_case}.  We have labeled the curves $\alpha_4$ and $\gamma$ which need to be used to complete the argument.

Assume that only one of the curves is separating.  Without loss of generality, we relabel the curves so that $\alpha_3$ bounds a subsurface which contains  $\alpha_1$ and $\alpha_2$.  If the subsurface containing these three curves is not a torus, then we can find curves $\alpha_4$ and $\gamma$ as above and apply the same argument.  If $\alpha_3$ bounds a torus containing the other two curves, then the result follows easily since $\alpha_1$ and $\alpha_2$ bound a cylinder of signed area $\Omega_1$, so
\[ [\alpha_1] + [\alpha_2] = \Omega_1 \rho \]  
while $\alpha_3$ bounds a torus whose area is exactly $\Omega - \Omega_1$.

The last remaining case is that each curve $\alpha_i$ bounds a surface $\Sigma_i$.  It is easy to decompose each $\Sigma_i$ into pairs of pants (each having exactly one separating boundary component), and inductively prove that
\[ [\alpha_i] = - \Omega_i \rho + \chi(\Sigma_i) T .\]
The result then easily follows from additivity of area and Euler characteristic, together with Lemma \ref{K_theory_relation} which we now prove.
\end{proof}

This is the only place where we have used immersed curves.  As we stated in the introduction, we could have avoided the use of immersed curves by finding the algebraic quasi-isomorphism which follows from the area-preserving isotopy which is implicit in comparing the final curves of Figure \ref{pair_pants_torus_proof}.  Since we did not include immersed curves in our Fukaya category, only the quasi-isomorphism which follows from that figure is in fact used in our arguments.

\begin{figure}
\begin{center}
$\begin{array}{c@{\hspace{.25in}}c} 
\input{genus_2_pair_of_pants.pstex_t} & \input{genus_2_test_curves.pstex_t}
\end{array}
$
\end{center}
\caption{}
\label{genus_2_case}
\end{figure}

\begin{lem} \label{K_theory_relation}
The subgroup of $K(\Fuk(\Sigma))$ generated by $\bR \rho$ and $T$ maps isomorphically to $\bR \oplus \bZ / \chi(\Sigma) \bZ $ under the holonomy and winding number maps.  In particular, 
\[ \Omega_{\Sigma} \rho - \chi(\Sigma) T = 0 .\]
\end{lem}
\begin{proof}
The winding number map is surjective since the image of a bounding curve under the winding number map determines the Euler characteristic of the surface that it bounds modulo $\chi(\Sigma)$, and we can certainly find unobstructed curves (or pairs of curves) bounding surfaces of Euler characteristic for an integer from $1$ to $\chi(\Sigma)-1$. To prove that the map $(\hol, \wind)$ is a surjection, we must find a twisted complex mapping to
\[ (\Omega, n) \]
for any $n \in \bZ/ \chi(\Sigma) \bZ$ and real number $\Omega$.  To achieve this, we simply choose a bounding curve $\gamma_0$ with winding number $n$ and add to it an element of the subgroup of $K$-theory generated by boundaries of cylinders as in Lemma \ref{cylinders}. We must now prove injectivity.

Assume that $\Sigma$ has even genus.  Then there exists a curve $\gamma$ which separates $\Sigma$ into two symplectomorphic surfaces.  The previous lemma (only the part which does not require the present lemma) implies that
\[ \gamma = \frac{\Omega_{\Sigma}}{2} \rho - \frac{\chi(\Sigma)}{2} T .\]
By considering a symplectomorphism which permutes the two surfaces, we conclude that $\gamma = \gamma[1]$, which implies that in the $K$-group,
\[ \Omega_{\Sigma} \rho - \chi(\Sigma) T = 0 .\]

In the odd genus case, instead of a separating curve, we consider a separating pair of curves to obtain the same relation.  Once we impose this relation, we obtain a group that is isomorphic to the desired product.
\end{proof}

We now understand the subgroup of $K$-theory generated by collections of curves which bound a subsurface.  We use this to prove the main proposition of this section.
\begin{proof}[Proof of Proposition \ref{homologous} ]
It suffices to express the different between the $K$-theory classes represented by any two homologous curves $\alpha_1$ and $\alpha_2$ as a linear combination of $T$ and $\rho$.  There is an element $\phi$ of the Hamiltonian Torelli group mapping $\alpha_1$ to $\alpha_2$.  Since this group is generated by Dehn twists about separating curves $\gamma_i$ and by pairs of twists about bounding curves $\beta_j$ and $\beta'_j$, Lemma \ref{dehn-ext} proves
\[ [\alpha_2] = [\phi(\alpha_1)] = [\alpha_1] + \sum_{j} [\beta_j] - [\beta'_j] .\]
By the previous lemma the sum in right hand-side vanishes since it can be expressed in terms of $\rho$ and $T$, and we know its holonomy and winding number both vanish since the holonomies and winding numbers of $\alpha_1$ and $\alpha_2$ are equal.
\end{proof}

\section{Computing the  $K$-theory} \label{computing}
We are left with some technical details that need to be cleaned up in order to complete the
\begin{proof}[Proof of Theorem \ref{main}]
We have constructed a map to the desired direct sum of groups which is clearly surjective.  To prove injectivity, it will suffice to prove that a twisted complex mapping to the trivial homology class of $\Sigma$ lies in the subgroup of $K$-theory generated by $\rho$ and $T$.   Since we already know that the classes representing separating curves lie in the subgroup generated by $T$ and $\rho$, we may assume that $[\Gamma]$ can be expressed as a sum of classes corresponding to non-separating curves.

We now use the Lickorish generators of the mapping class group introduced in \cite{lickorish} (see Figure \ref{lickorish-gen}).  Since the mapping class group acts transitively on isotopy classes of non-separating curves, Lemma \ref{dehn-ext} implies that that the non-separating curves that appear in $[\Gamma]$ can be written as a linear combination of the Lickorish generators and their translates.  Note that the curves $\alpha_i$ and $\beta_i$ are a basis for $H_1(\Sigma, \bZ)$, and $\gamma_i$ bounds a pair of pants with two of these basis elements, so that
\[[\alpha_i] - [\alpha_{i+1}] + [\gamma_i] = T + \Omega_i \rho \]
which $\Omega_i$ is the area of the pair of pants.  In particular, if $[\Gamma]$ represents the trivial class in homology, then we can decompose it as a linear combination of the left hand side of the above equation (these relations suffice), which implies that $\Gamma$ does indeed lie in the subgroup generated by $T$ and $\rho$.

\end{proof}

\begin{figure}
\begin{center}
\input{lickorish_generators.pstex_t}
\end{center}
\caption{}
\label{lickorish-gen}
\end{figure}

\appendix
\section{Winding number and Holonomy}\label{wind_holo}
Let $\gamma$ be an immersed piecewise smooth curve in a surface $\Sigma$.  We may therefore write $\gamma$ as a union of smooth segments 
\[I_{m}:[0,1] \to \Sigma \]
which intersect at corners $c_m$ lying on $I_{m} \cap I_{m+1}$.  Let $\theta_m$ be the angle from the tangent vector of $I_m$ to the tangent vector of $I_{m+1}$.  

We can trivialize the restriction of the tangent space of $\Sigma$ to each segment $I_m$ with the tangent vector of $I_m$ corresponding to the $x$-axis.  A smooth vector field $\nu$ which does not vanish on $I_m$ therefore yields a map from the closed interval to $S^1$.  This map lifts to a map $\tilde{\nu}_{m}: I \to \bR$.  
\begin{defin}
The {\bf winding number} of $\nu$ around the curve $\gamma$ is the integer
\[ \wind_{\nu}(\gamma) = \sum_{m}\tilde{\nu}_{m}(1) - \tilde{\nu}_{m}(0) - \theta_{m} .\]
\end{defin}

We now specialize to the case that $\nu$ has a unique zero on $\Sigma$.  It is well known that the degree of vanishing of $\nu$ must equal the Euler characteristic of $\Sigma$.  In particular, if we choose a small disc centered at the vanishing point, and orient its boundary circle $\gamma$ as the boundary of the complement, then we can explicitly check that
\[ \wind_{\nu}(\gamma) = 1 - \chi(\Sigma) .\]

The following results are classical, and detailed discussions can be found in \cite{reinhart1}, \cite{reinhart2}, \cite{chillingworth1}, and \cite{chillingworth2}.
\begin{thm}
\mbox{}
\begin{itemize}
\item Modulo the Euler characteristic, the winding number is an invariant of the isotopy class of a curve.
\item If $\Sigma' \subset \Sigma$ is an oriented subsurface which does not include the zero of $\nu$ and which has boundary $\cup_{i} \gamma_i$, then
\[ \chi(\Sigma') = - \sum_{i} \wind_{\nu}(\gamma_i) \]
\end{itemize}
\end{thm}
\begin{proof}[Sketch of proof]
The first result can be checked locally by observing that the winding number changes by $\pm \chi(\Sigma)$ when a family of curves ``passes through the zero.''  To prove the second part, we fill in the boundary components of $\Sigma'$ with discs, and extend $\nu$ to a vector field $\bar{\nu}$ on this closed surface $\overline{\Sigma}'$.  The winding number of $\nu$ around each boundary component determines the (weighted) number of zeros in $\overline{\Sigma}'$ allowing us to compute the Euler characteristic of this surface in terms of the winding numbers around $\gamma_i$.  On the other hand, the Euler characteristic of $\Sigma$ and its number of boundary components also determine the Euler characteristic of $\bar{\Sigma'}$. 
\end{proof}

The other auxiliary choice we make is that of a $1$-form on the unit tangent bundle of $\Sigma$ whose differential equals to pull-back of $\omega$.  Such a $1$-form corresponds to a connection $A$ on the tangent bundle whose curvature gives our desired symplectic form.  This choice constrains the area of $\Sigma$, but it is obvious that the Fukaya category is insensitive to the total area.  Given a segment $I$ of a piecewise smooth curve $\gamma$ the tangent vector and the connection $A$ induce two trivializations of the tangent space.  Comparing them allows us to define a real number
\[ \hol_{A}(I) .\]
As above, taking the sum over all segments $I$ and adding the angles at every vertex allows us to define the holonomy of any piecewise smooth curve $\gamma$.  If $\gamma$ were smooth, the above holonomy would equal the integral of the connection $1$-form over the unique lift of $\gamma$ to the unit sphere bundle.
\begin{lem}
If $\gamma$ and $\gamma'$ bound a cylinder of area $\Omega$, then
\[\hol_{A}(\gamma) - \hol_{A}(\gamma') = \Omega. \]
\end{lem}
\begin{proof}
This is an immediate consequence of Stokes's theorem.
\end{proof}

\section{Symplectomorphism Groups}\label{symplecto_gps}
We collect some elementary facts about the group of symplectomorphisms of a surface and its various sub-quotients.  For an overview of these facts from a slightly different point of view, see \cite{KM}.

\begin{defin}
The {\bf Hamiltonian mapping class group} of a symplectic manifold $\Sigma$ is the quotient
\[ \Symp(\Sigma)/\Ham(\Sigma). \]
\end{defin}

We will now construct a set of generators for this group.  Let $\gamma$ be a simple closed curve in $\Sigma$.  Choose a symplectomorphism between a collar neighbourhood of $\gamma$ (where the two sides of the collar have equal area) and an annulus of area $2\epsilon$ such that the image of $\gamma$ is the core curve. Choose a $\delta < \frac{\epsilon}{2}$.  Let us give the annulus radial coordinate $r \in [-\epsilon, \epsilon]$ and angle coordinate $\theta \in [0,1]/\{0,1 \}$ and choose a $C^{\infty}$ function $f$ on $\bR$ with even derivative such that
\begin{align*}
f(r) = 0 & \iff r \leq - \delta
\\ f(r) = 1 & \iff r \geq \delta.
\end{align*}
\begin{defin}
The map
\[ (r,\theta) \mapsto ( r , \theta + f(r) ) \]
is a symplectic Dehn twist about $\gamma$.
\end{defin}
Note that the symplectic Dehn twist is constant away from the width-$\delta$ collar of $\gamma$, hence extends to $\Sigma$.  

\begin{lem}
Any two symplectic Dehn twists about $\gamma$ are Hamiltonian isotopic
\end{lem}
\begin{proof}[Sketch of Proof]
Any two collar neighbourhoods of the same area are Hamiltonian isotopic, so one can prove the result for the annulus by explicitly producing the desired Hamiltonian function.  Alternatively, one can use the characterization of the group of Hamiltonian symplectomorphisms as the Kernel of the Flux homomorphism, and check that the flux of a Dehn twist followed by the inverse of a Dehn twist vanishes.\end{proof}

We therefore have a well defined lift of the Dehn twist about $\gamma$ to the Hamiltonian mapping class group.  We will use $\tau_{\gamma}$ for this symplectomorphism.
\begin{lem} \label{symp-dehn-gen}
The Hamiltonian mapping class group is generated by symplectic Dehn twists.
\end{lem}
\begin{proof}
Moser's theorem implies that we have a short exact sequence,
\[ 1 \to \Symp_0(\Sigma) \to \Symp(\Sigma) \to Mod(\Sigma) \to 1 .\]
Taking the quotient of the first two groups by $\Ham(\Sigma)$, and using the Flux homomorphism to identify the first quotient with first cohomology, we obtain:
\[ 0 \to H^1(\Sigma, \bR) \to \Symp(\Sigma)/\Ham(\Sigma) \to Mod(\Sigma) \to 1 .\]
Since Dehn twists generate the mapping class group, it suffices to show that the subgroup of $ \Symp(\Sigma)/\Ham(\Sigma) $ generated by symplectic Dehn twists contains the first term in the above exact sequence.  Indeed, if $\gamma$ and $\gamma'$ are non-separating simple closed curves which are isotopic  through a cylinder of area $\Omega$, then
\[ \phi = \tau_{\gamma} \circ \tau^{-1}_{\gamma'} \]
is an element of $ \Symp(\Sigma)/\Ham(\Sigma) $ which represents the Poincar\'e dual of $\Omega \times [\gamma]$ in $H^1(\Sigma, \bR)$.  Such elements clearly generate $H^1(\Sigma, \bR)$.
\end{proof}
\begin{rem}
Note that is suffices to have lifts of Dehn twist about curves which generate the mapping class group.  In particular, the lifts of curves which are isotopic to the Lickorish elements generate the Hamiltonian mapping class group.
\end{rem}

\begin{defin}
The {\bf Hamiltonian Torelli group} is the subgroup of the Hamiltonian mapping class group which acts trivially on the homology of $\Sigma$.
\end{defin}

This group also has a convenient set of generators.
\begin{cor} \label{tor-twist-gen}
The Hamiltonian Torelli group is generated by symplectic Dehn twists about bounding curves and by bounding symplectic twists.
\end{cor}
\begin{rem}
Recall that a bounding twist is an element of the mapping class group which is the composite of a Dehn twist along a curve $\alpha_1$ with the inverse Dehn twist along a curve $\alpha_2$ where $\alpha_1$ and $\alpha_2$ form the boundary of an embedded submanifold of $\Sigma$.
\end{rem}
\begin{proof}
This is an immediate consequence of the proof of the previous Lemma, and a classical result of Powell, \cite{powell}, that the Torelli group is generated by Dehn twist along separating curves and by bounding twists.
\end{proof}

\section{$A_{\infty}$ pre-categories} \label{precategories}
Our definition of $A_{\infty}$ category differs from the usual notion since we do not work with $\bZ$-graded complexes.  To fix the notation, we will say that an $A_{\infty}$ category $\cA$ consists of a set of objects $\Ob$ and morphisms $\Mor(A,B)$ between any two objects which are $\bZ / 2 \bZ$ graded vector spaces, together with a collection of maps
\[ m_{n} : \Mor(A_{n-1}, A_{n}) \otimes \cdots \otimes \Mor(A_{0}, A_1) \to \Mor(A_0, A_n) \]
satisfying the $A_{\infty}$ equation which appears in Lemma \ref{infty_relation}.  For more on $A_{\infty}$-categories, see \cite{lefevre}.

The $A_{\infty}$ categories that are usually studied satisfy an additional axiom analogous to the existence of a unit morphism in an ordinary category.  However, the Fukaya category of a surface is not an $A_{\infty}$ category.  In particular, since a curve is never transverse to itself, the morphism space $CF^*(\gamma, \gamma)$ is never defined, so it doesn't make sense to discuss the existence of units from a naive point of view.  We will adopt the framework of $A_\infty$-pre-categories introduced by Kontsevich and Soibelman \cite{KS1}.  

\begin{defin}
An $A_{\infty}$ pre-category consists of the following data:
\begin{itemize}
\item A set of objects $\Ob$.
\item For each integer $ k \geq 1$, a set of transversal objects $\Ob_{k}^{tr} \subset \Ob^{k}$ satisfying the recursive conditions that under the projection map (which forgets any $k-i$ elements),
\[ Im(\Ob^{tr}_{k}) \subset \Ob^{tr}_{i} .\]
\item For each transversal pair $(A,B)$, a $\bZ/ 2 \bZ$ graded vector space $\Mor(A,B)$.
\item A collection of maps $m_n$ for each transversal multiplet $(A_0, \cdots ,A_n)$
 \[ m_{n} : \Mor(A_{n-1}, A_{n}) \otimes \cdots \otimes \Mor(A_{0}, A_1) \to \Mor(A_0, A_n) \]
which satisfy the $A_{\infty}$ equation.
\end{itemize} 
\end{defin}

Note that the last condition only makes sense because we've required the set of transversal objects to satisfy the appropriate recursive property.  We have the following notion of quasi-isomorphism.
\begin{defin}
An element $f$ of $\Mor(A,A')$ satisfying $m_1(f)=0$ is called a {\bf quasi-isomorphism} if the maps
\begin{align*}
 m_2(\_ , f) :& \Mor(A', B') \to \Mor(A, B') \\
 m_2(f, \_):& \Mor(B,A)  \to \Mor(B,A')
\end{align*}
are quasi-isomorphisms of complexes (i.e. induce isomorphisms on homology) whenever $(B,A,A')$ or $(A,A',B')$ are transverse triples.
\end{defin}
An $A_{\infty}$ pre-category is said to be unital if there are enough quasi-isomorphisms.  Formally, we require
\begin{defin}
An $A_{\infty}$ pre-category is {\bf unital} if for any object $A$, and a finite collection of transversal sequences $\{ S_i \}$ there exists objects $A^+$ and $A^-$, quasi-isomorphic to $A$ such that the sequences
\[ \{ A^- , S_i , A^+ \} \]
are transversal for every $i$.
\end{defin}

We note that in a unital $A_{\infty}$ pre-category, quasi-isomorphism is an equivalence relation among transversal objects.
\begin{lem}
Assume that $\{ A_1, A_2, A_3 \}$ are a transversal triple such that 
\begin{align*}
e_{1,2} & \in \Mor(A_1, A_2) \\
e_{2,3} & \in \Mor(A_2, A_3)
\end{align*}
are quasi-isomorphisms.  Then
\[ e_{1,3} = m_2(e_{2,3}, e_{1,2}) \in \Mor(A_1, A_3) \]
is also a quasi-isomorphism.
\end{lem}
\begin{proof}
The fact that $m_2$ is a chain map implies that $e_{1,3}$ is a closed morphism.  Choose an arbitrary $B$ which is transversal to $(A_1, A_3)$.  It may not be transverse to $A_2$.  However, there exists $B'$ with a quasi-isomorphism $e \in \Mor(B',B)$ such that the sequences $(B' , B ,A_1, A_3)$ and $(B', A_1, A_2, A_3)$ are transverse.  The result now follows from the fact that the following diagram commutes up to homotopy, and all but one of its arrows are known to induce isomorphisms on homology
\[
\xymatrix{  \Mor(B',A_1) \ar[rr]^-{m_2(e_{1,2}, \_)} & & \Mor(B', A_2) \ar[rr]^-{m_2(e_{2,3}, \_)} & & \Mor(B',A_3)
\\   \Mor(B,A_1) \ar[rrrr]^-{m_2(e_{1,3}, \_)} \ar[u]^-{m_2(\_, e)} & & & & \Mor(B,A_3). \ar[u]^-{m_2(\_, e)}}
\] 
\end{proof}

\section{Twisted complexes and $K$-theory} \label{twisted-and-K}
We will use twisted complexes to derive any $A_{\infty}$ pre-category.  Since this is the same construction as for $A_{\infty}$ categories, we follow \cite{fukaya3} closely, and the only concern is to define transversal sequences.

We begin by constructing a category which is closed under shifts. The shift $A[1]$ of an object $A$ has morphisms defined by
\[ \Mor(A[1],B) = \Mor(A,B)[1] \]
where the right hand side is the $\bZ/2\bZ$-graded vector space whose graded pieces are exactly those of $\Mor(A,B)$ shifted by $1$, and similarly for morphisms with target $A[1]$.  A sequence containing $A[1]$ is transversal if and only if the corresponding sequence with $A$ is transversal.

From now on, we will assume that the pre-category $\cA$ is closed under shifts, since we can enlarge it using the previous construction.
\begin{defin}
A {\bf twisted complex} $C$ over an $A_{\infty}$-pre-category $\cA$ is a finite set of objects  $ \{ A_i \}_{i=0}^{n}$ together with a collection morphisms of degree $1$ which we write as $\Delta =( \delta_{i,j})_{i<j}$ such that
\[ \sum_{i_1 < \cdots < i_{k}}  m_{k}(\delta_{i_{k-1}, i_k}, \cdots, \delta_{i_{2}, i_1}) =0. \] 
\end{defin}

Given twisted complexes $A = (A_i, \Delta^{A})$ and $B = (B_k, \Delta^{B})$, such that every sequence of objects
\[( A_{i_1}, \cdots, A_{i_d}, B_{k_1}, \cdots  ,B_{k_e}) \]
is transverse, we define
\[ \Mor(A,B) = \bigoplus_{i,k} \Mor(A_i, B_k) .\]
The morphisms $\Delta^{A}$ and $\Delta^{B}$ induce a differential on this graded module.  Given an element $F = (f_{i,k})$ of the above module we define
\[ \hat{m}_{1}(F)_{i,k} = \sum_{i<i_1 < \cdots< i_{d-1} < i_d, k_1 < \cdots <  k_e <k} m_{d+e} (\delta^{B}_{k_{e}, k}, \cdots , \delta^{B}_{k_1,k_2}, f_{i_d,k_1}, \delta^{A}_{i_{d-1}, i_d}, \cdots, \delta^{A}_{i, i_1}),\]
and extend $\hat{m}_{1}$ linearly.  As a shorthand, we will write $I$ for an increasing sequence $i_1 < \cdots< i_d < i_{d+1}$, and $\Delta^{A}_{I}$ for the corresponding sequence $(\delta^{A}_{i_d,i_{d+1}}, \ldots, \delta^{A}_{i_1, i_2})$.

\begin{defin}
Given an increasing sequence of integers $I$, its {\bf length} (written $|I|$) is the number of elements.
\end{defin}

Using this notation, the above formula becomes
\[ \hat{m}_1(F) = \sum_{I, K} m_{|I|+|K|+1}(\Delta^{B}_K, f_{i,k}, \Delta^{A}_I) \]
where the sum is taken over all meaningful choices of sequences $I$ and $K$; in particular, the last morphism in $I$ has $A_i$ as target, while the first morphism in $K$ has $B_k$ as source.  We refer to \cite{fukaya3} for the proof that $\hat{m}_1$ squares to zero and for the proof of the $A_{\infty}$ category analogue of
\begin{lem}[Definition-Lemma]
If $\cA$ is an $A_{\infty}$ pre-category, then twisted complexes over $\cA$ also form an $A_{\infty}$ pre-category $Tw(\cA)$ satisfying the following conditions:
\begin{itemize}
\item Objects are twisted complexes with transversal pairs and morphisms between them defined as above.
\item A sequence $(C_{i} = (c_{i,j}))_{i=1}^{n}$ of twisted complexes is transverse if and only if every sequence of the form
\[ (c_{1,1}, \cdots, c_{1,j_1},\cdots, c_{i,1}, \cdots, c_{i,j_i}, \cdots, c_{n,1}, \cdots, c_{n,j_n}) \]
is transverse.  It is understood that the second index in each subsequence $( c_{i,1}, \cdots, c_{i,j_i})$ is always increasing.
\item The higher products $\hat{m}_k$ are defined among transversal objects in analogy with $\hat{m}_1$ (We omit the explicit formulae which appear in \cite{fukaya3}).
\end{itemize}
Further, if $\cA$ is unital, then so is $Tw(\cA)$. \noproof
\end{lem}

The only higher product for which we will need a formula is $\hat{m}_2$.  If $(A,\Delta^{A})$, $(B, \Delta^{B})$ and $(C, \Delta^{C})$ is a transverse triple of twisted complexes, then
\[ \hat{m}_{2}(G,F) = \sum_{N,K,I} m_{|N|+|K|+|I|+2} ( \Delta^C_{N}, g_{l,n}, \Delta^{B}_{K}, f_{i,k}, \Delta^{A}_{I} )  \]
for $F=(f_{i,k}) \in \Mor(A,B)$ and $G=(g_{l,n}) \in \Mor(B,C)$, and where the sum is taken, as usual, over all possible increasing sequences of labels $N$, $K$ and $I$.

The main advantage of the category of twisted complexes is that it is triangulated.  This is proved in \cite{fukaya3} for example.  We will restrict ourselves to explaining what the distinguished triangles are, and will not prove all the analogous results in the pre-category setting.  

Given a closed morphism of twisted complexes $F \in \Mor(A,B)$, we construct $Cone(F)$ as the twisted complex consisting of the sequence of objects
\[ (A_1, \cdots, A_n, B_1, \cdots, B_m) \]
and with morphisms
\[ \delta_{i,j} = \begin{cases} \delta^{A}_{i,j} & \mbox{if } i,j \leq n \\
f_{i,j-n} & \mbox{if } i \leq n \mbox{ and }  j > n
\\ \delta^B_{i-n,j-n} & \mbox{if } i , j > n
\end{cases}
.\]

\begin{defin}
Let $\cA$ be an $A_{\infty}$-pre-category, $G$ the free abelian group generated by quasi-isomorphism classes of objects in $Tw(\cA)$, and $H$ the subgroup of $G$ generated by elements $[A]+[B]-[C]$ where $C$ is the cone of a closed morphism from $A$ to $B$.  The $K$-theory of $\cA$ is the quotient
\[ K(Tw(\cA)) = G/H .\]
\end{defin}

This ad-hoc definition can be formalized by introducing the appropriate notion of triangulation on an $A_{\infty}$ pre-category.  We will not take up this task here since it would involve rather technical homological algebra.  We complete this appendix by proving some relatively straight forward results that will allow us to simplify the computation of $K(\Fuk(\Sigma))$.

\begin{lem}
$K(\cA)$ is generated by objects of $\cA$.
\end{lem}
\begin{proof}
Every twisted complex is an iterated cone on objects of $\cA$.
\end{proof}
\begin{lem}
Assume that $I: A \to A'$ is a closed morphism in $\cA$ such that the maps
\begin{align*}
\hat{m}_2( \_ , I): \Mor(C,A) &\to \Mor(C, A') \\
\hat{m}_2( \_ , I): \Mor(D,A) & \to \Mor(D,A') \end{align*}
are both isomorphisms in homology.  If $F: C \to D$ is any closed morphism, then 
\[ \hat{m}_2( \_ , I): \Mor(Cone(F),A) \to  \Mor(Cone(F), A') \]
is an isomorphism in homology.
\end{lem}
\begin{proof}
We apply the $5$ lemma to the diagram:
\[
\xymatrix{  H_{i-1}(D,A) \ar[r] \ar[d] & H_{i}(C,A) \ar[r] \ar[d] & H_{i}(Cone(F),A)  \ar[r] \ar[d] &  H_{i}(D,A)  \ar[r] \ar[d] &   H_{i+1}(C,A) \ar[d]  \\
 H_{i-1}(D,A') \ar[r] & H_{i}(C,A') \ar[r]  & H_{i}(Cone(F),A')  \ar[r] &  H_{i}(D,A')  \ar[r] &   H_{i+1}(C,A')   \\
}
\] 
\end{proof}
\begin{cor}
If $A$ and $A'$ are quasi-isomorphic in $\cA$, then they are quasi-isomorphic in $Tw(\cA)$.
\end{cor}
\begin{cor}
Let $F: A \to A'$ be a morphism between objects of $Tw(\cA)$ such that $m_2( \_ , F)$ and $m_2(F, \_)$ induce chain isomorphisms between $\Mor(B,A)$ and $\Mor(B,A')$ and between $\Mor(A',C)$ and $\Mor(A,C)$ whenever $A$ and $B$ are objects of $\cA$ for which the triples $\{ B, A, A' \}$ and $\{ A , A' , C \}$ are transverse.  Then $A$ and $A'$ are quasi-isomorphic in $Tw(\cA)$.
\end{cor}
\begin{rem}
This means that it suffices to test a quasi-isomorphism using objects of $\cA$, instead of testing it on arbitrary twisted complexes.
\end{rem}

\end{document}